\documentclass[lettersize,journal]{IEEEtran}
\usepackage{amsmath,amsfonts,amsthm, amssymb}
\usepackage{mathtools}
\usepackage{algorithm}
\usepackage{array}
\usepackage[caption=false,font=normalsize,labelfont=sf,textfont=sf]{subfig}
\usepackage{textcomp}
\usepackage{stfloats}
\usepackage{url}
\usepackage{verbatim}
\usepackage{graphicx}
\usepackage{cite}
\usepackage{amssymb}
\usepackage[figuresright]{rotating}
\usepackage{booktabs}
\usepackage{algorithmicx}
\usepackage{algpseudocode} 
\usepackage{booktabs}
\usepackage{url}
\usepackage{multirow}
\usepackage[table,xcdraw]{xcolor}
\usepackage{color}
\usepackage{soul}

\definecolor{mistyrose}{rgb}{1.0, 0.89, 0.88}
\newcommand{\net}{NH-CSR}
\newcommand{\dx}{\,\mathrm{d}x}
\newcommand{\R}{\mathbb{R}}
\newcommand{\calT}{\mathcal{T}}

\begin{document}

\title{Multiscale Corrections
by Continuous Super-Resolution}

\author{Zhi-Song~Liu, 
        Roland Maier,
        and~Andreas~Rupp
\thanks{Zhi-Song Liu is with Lappeenranta--Lahti University of Technology (LUT), Finland. Email: zhisong.liu@lut.fi}
\thanks{Roland Maier is affiliated with the Institute for Applied and Numerical Mathematics at Karlsruhe Institute of Technology (KIT), Germany}
\thanks{Andreas Rupp is affiliated with the Department of Mathematics at Saarland University, Germany. Email: andreas.rupp@uni-saarland.de}}

\markboth{}
{Shell \MakeLowercase{\textit{et al.}}: A Sample Article Using IEEEtran.cls for IEEE Journals}

\IEEEpubid{}

\maketitle

\begin{abstract}
 Finite element methods typically require a high resolution to satisfactorily approximate micro and even macro patterns of an underlying physical model. This issue can be circumvented by appropriate 
 multiscale strategies that are able to obtain reasonable approximations on under-resolved scales. In this paper, we study the implicit neural representation and propose a continuous super-resolution network as a
 correction strategy for multiscale effects. It can take coarse finite element data to learn both in-distribution and out-of-distribution high-resolution finite element predictions. Our highlight is the design of a local implicit transformer, which is able to learn multiscale features. We also propose Gabor wavelet-based coordinate encodings which can overcome the bias of neural networks learning low-frequency features. Finally, perception is often preferred over distortion so scientists can recognize the visual pattern for further investigation. However, implicit neural representation is known for its lack of local pattern supervision. We propose to use stochastic cosine similarities to compare the local feature differences between prediction and ground truth. It shows better performance on structural alignments. Our experiments show that our proposed strategy achieves superior performance as an in-distribution and out-of-distribution super-resolution strategy.
\end{abstract}


\section{Introduction}
\label{Introduction}

\noindent This study paves a new direction of using deep learning-based continuous super-resolution as a numerical homogenization tool. To this end, we address the issue of efficiently simulating physical phenomena characterized by effects on several length scales.
In this context, the term \emph{numerical homogenization} refers to a general mathematical strategy that allows to employ fine-scale features of a solution (sometimes referred to as oscillations) in a coarse-scale simulations. Thus, numerical homogenization can in some sense be translated to \emph{multiscale correction} in what follows. A comparison between \emph{classical} mathematical homogenization approaches and our super-resolution-based approach is illuminated in the below paragraphs.

\noindent We denote our approach \emph{continuous}, since as opposed to most super-resolution approaches, which are trained on fixed integer scales (e.g., the super-resolution is conducted by factors of 2, 4, \dots), we allow for an arbitrary scaling and also allow for out-of-distribution resolutions. That is, we may train our model to downscale the simulation \( 16 \) times, and use it to downscale our numerical results \( 10 \) or \( 20 \) times.

Examples of multiscale phenomena with effects on multiple scales are flow and transport through porous media. To reliably predict such phenomena, information on the micro-scale is necessary. However, if, e.g., flow through soil is considered, we are usually only interested in some macroscopic outcome, such as the final macroscopic location of some contaminant in a field. In mathematical terms, such processes are described by partial differential equations (PDEs). More precisely, the equilibrium state of a flow process is characterized by an unknown function $u$ that solves the PDE
\begin{equation}\label{eq:model}
 \left\{
 \begin{aligned}
  -\operatorname{div} (A\nabla u) &= f&& \quad\text{ in } \Omega,\\
  u &= 0&& \quad\text{ on } \partial \Omega,
 \end{aligned}
 \quad\right.
\end{equation}
in a bounded domain $\Omega \subset \R^d$. Here, $A \colon \Omega \to \R$  is a scalar-valued coefficient, and $f \colon \Omega \to \R$ is a source or sink. The multiscale character is encoded in the structure of $A$, which is strongly heterogeneous on a fine scale.

\noindent Our main assumption in this manuscript is that \( \Omega \) is a rectangular domain covered by some (potentially very fine) pixel mesh on which \( A \) varies between two discrete values. The restriction of \( A \) can be justified in many cases, e.g., if flow through a porous medium, which results in to two distinct diffusivities, and by a measurement device that is only capable of classifying the medium on the pixel mesh.

\noindent Since problems of the type \eqref{eq:model} are not only relevant in soil science~\cite{ZHANG2023105247,Rupp19}, oil recovery~\cite{RainerOil} (and the references therein), drug development~\cite{https://doi.org/10.1002/zamm.201200196}, and medicine (where the porous medium could, e.g., be a brain~\cite{https://doi.org/10.1002/pamm.200810201} or the skin of an individual~\cite{Reis_2015}), but in many disciplines of science and engineering, many approaches in the literature try to tackle them.

\noindent \textbf{Mathematical approaches.} From a mathematical perspective, \eqref{eq:model} cannot be approximated by standard numerical schemes such as finite element, finite volume, or finite difference schemes on the coarse scale of interest since their respective meshes would need to resolve the fine-scale structure. This would result in an unfeasible number of degrees of freedom and computational demands that forbid the application of these techniques to technically relevant problems. However, there are mathematically rigorous and efficient approaches to tackle this issue, such as~\cite{HouW97,MatS02,EE03,EE05}. These approaches are based on analytical homogenization theory, see~\cite{BlaLB23,Tar09}, and the references therein, that require the coefficient to exhibit specific properties (such as periodicity or scale separation) that are rarely met in practice.  Alternatively, heuristic approaches try to infer easy-to-compute homogeneous problems from statistical considerations \cite{RayRSK18}.

\IEEEpubidadjcol

\noindent A promising and reliable alternative to the approaches above are computational multiscale methods, referred to as \emph{numerical homogenization}. The idea is to generate optimally adapted finite element solutions from a coarse-scale finite element solution enriched by suitable fine-scale information. In particular, such approaches provably work for a vast class of heterogeneous media with minimal structural assumptions, see, e.g., \cite{EfeGH13,ChuEL18b,BabL11,MaSD22,MaS22,HenP13,MalP14,MalP20,Owh17,OwhS19,Mai21,DonHM23,AltHP21}. Although these methods are highly efficient for solving a multiscale problem once the enrichment functions are computed, they need substantial pre-computations per choice of $A$.

\noindent\textbf{Our approach.} We propose a different approach in the spirit of numerical homogenization. The goal is to learn a mapping from a (potentially incorrect) coarse-scale solution to an improved (upscaled) solution that appropriately represents fine-scale features for a large class of coefficients $A$ and not for each $A$ individually as in the above techniques. That is, we connect the PDE problem to image super-resolution, which is a fundamental task in image processing to recover fine textures in given images. As tools for ill-posed inverse problem, deep learning-based super-resolution approaches have shown remarkable progress in the past decade and have been widely used in digital entertainment, content creation, and fashion design. Recently, researchers have been interested in applying super-resolution to physical problems, including fine-grained climate prediction~\cite{climate}, remote sensing~\cite{remote}, and coastal simulations \cite{LiuBAR24}. Our contributions can be summarized as follows:
\begin{itemize}
 \item We propose a numerical homogenization strategy via continuous super-resolution (\net), incorporating a coarse (and possibly) inaccurate approximation of~\eqref{eq:model} and the coefficient map $A$ to achieve reliable approximations via super-resolution.
 \item Our proposed framework achieves better visual quality due to wavelet implicit neural representations (WIRE), which use a continuous real Gabor wavelet activation functions as position-based image functions. It better concentrates on space frequency and has excellent biases for representing images.
 \item We also employ the multiscale implicit image function (MS-IIF) for frequency interpolation. It can faithfully preserve low-frequency information while introducing high-frequency details for better visualization. 
 \item To encourage the neural network to learn spatial structural information, we introduce a novel stochastic cosine similarity loss that can grab the non-local pixel correlation in the 2D space, which gains significant improvements in super-resolution.
 \item We empirically demonstrate the supremacy of the proposed \net\ over the existing state-of-the-art techniques in quantitative and qualitative measurements. We further experiment with ultra-resolution schemes to show that our method is superior to others in out-of-distribution super-resolution.
\end{itemize}

\begin{figure*}[t!]\centering
 \centerline{\includegraphics[width=\textwidth]{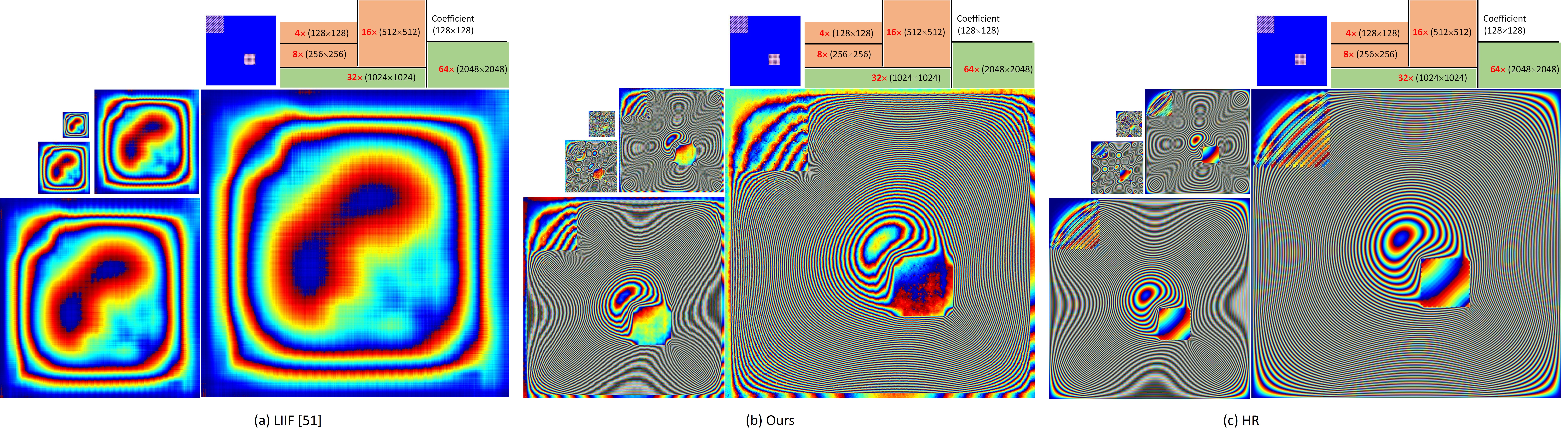}}
 \caption{\small \textbf{Continuous super-resolution for finite element data.} We use the proposed method to apply multiscale upscaling to the low-resolution image, both in-distribution (orange region, upscaling factor not larger than 16) and out-of-distribution (green region, upscaling factor larger than 16). We see from the comparison that ours achieves better visualization than the state-of-the-art LIIF~\cite{LIIF} approach. Note that the actual values are first normalized to [0, 255], then multiplied by a factor of 64, and taken modulo 255 to again be in the value range between 0 and 255. Finally, the result is converted to the RGB color space using the Matplotlib colormap to better visualize differences. This procedure aims at highlighting fine differences in simulation results and leads to the observed wavy patterns.}
 \label{fig:teaser}
\end{figure*}

\noindent Figure~\ref{fig:teaser} shows an exemplary comparison of our approach with local implicit image functions (LIIF)~\cite{LIIF} on the finite element data super-resolution. We apply the \net\ to the low-resolution data and produce multiscale super-resolved images. Since the model is trained for the setting with an upsampling factor of 16, any upsampling factors that are larger than 16 prescribe an out-of-distribution super-resolution. In the right upper corner (left), we also show the coefficients that are used for the data generation. The remaining plots show visualizations of the upsampled solutions, which are normalized to $[0,255]$, then multiplied by a factor 64 and cut again to the range $[0,255]$ by taking the values modulo 255. These are then used for a coloring using the RGB color space. This is done to highlight pattern differences more clearly. In particular, this leads to the oscillatory structure of the pictures. Overall, we observe that our approach outperforms LIIF in all upsampling scenarios.

\noindent \textbf{Theoretical grounding.} The computational multiscale methods, mentioned above, use the fact that the fine-scale behavior of the solution \( u \) can be reconstructed from  the coefficient\( A \). In particular, they correct coarse-scale shape functions in an \( A \)-dependent fashion to account for the fine-scale behavior. Afterward, they solve a coarse finite element problem with these `corrected' ansatz functions.

\noindent Contrasting that approach, we first solve the coarse finite element problem with `uncorrected' shape functions, and correct the solution \( u_h \) afterward. To that end, our approach receives the coefficient \( A \) and the approximated solution \( u_h \) to generate the fine-scale features. It can, thus, be justified similarly to its mathematical analogs.

\noindent In essence, since the PDE of interest is known, only the coefficient $A$ (and right-hand side \( f \)) is required to obtain a reasonable approximation, see also classical theory on finite element methods, e.g., \cite{Braess,BrennerScott,Bartels,KnabnerA21}. In a similar manner, several works have shown (theoretically and experimentally) that the coefficient-to-solution map can be well-approximated by neural networks, see, e.g., \cite{GeiPRSK21,BhaHKS21,GaoSW21,KutPRS19}. In that sense, our approach even works with more data than theoretically required, as the (easy-to-compute) coarse-scale solution is also taken as input, giving already a very rough approximation of the desired solution as a starting point.

\noindent \textbf{Structure.} The manuscript is structured as follows: In the next section, we discuss the available mathematical tools to tackle~\eqref{eq:model}. Afterward, we mention related works in the context of super-resolution. Section~\ref{Approach} details our approach regarding the possible application of super-resolution to achieve similar results as numerical homogenization approaches. Finally, Section~\ref{Experiments} comprises illustrations of our numerical results.

\section{Background and related  Works}
\label{related_works}

\subsection{Mathematical approaches}

The basis of numerical homogenization strategies are traditional discretization methods, such as the finite element method. We briefly state the definition of the method and then review some multiscale methods based on it. 

\noindent \textbf{Finite element approximations.} 
A classical way to approximate solutions to~\eqref{eq:model} is the finite element method. It is based on a weak formulation: multiplying~\eqref{eq:model} with a so-called test function and using integration by parts results in seeking a function $u$ that solves
\begin{equation}\label{eq:modelweak}
 a(u,v) : = \int_\Omega A \nabla u \cdot \nabla v \dx = \int_\Omega f v \dx =: (f,v)
\end{equation}
for all appropriate test functions $v$. The finite element method restricts possible functions to approximate $u$ (trial functions) and functions $v$ (test functions) to a finite-dimensional subspace that is built from a decomposition of the domain into a grid $\calT_h$ consisting of small boxes $K$ with side length $h > 0$. More precisely, let
\begin{equation*}
 V_h := \left\{\begin{aligned} &v\colon \Omega \to \R \;\colon\; v\vert_K \text{ is a polynomial of}\\& \text{partial degree } \leq 1 \text{ for all }K \in \calT_h\end{aligned}\right\} \cap C^0(\Omega), 
\end{equation*}
which is the space of continuous functions that are polynomials of degree at most one in each coordinate direction when restricted to a box of the grid $\calT_h$. The finite element method then seeks a function $u_h \in V_h$ such that
\begin{equation}\label{eq:fem}
 a(u_h,v_h) = (f,v_h)
\end{equation}
for all functions $v_h \in V_h$. Since $V_h$ is a finite-dimensional space, \eqref{eq:fem} is equivalent to solving a system of linear equations.

\noindent As mentioned above, the finite element method only provides reasonable approximations if $h < \varepsilon$, where $\varepsilon$ is the most miniature scale on which the coefficient varies. However, resolving all details can be a computationally demanding task or even not possible at all for large systems of equations. Therefore, alternative strategies are necessary.necessary

\noindent \textbf{Multiscale strategies.} 
An alternative to globally resolving fine scales are homogenization-type approaches. Suppose $A$ has a specific structure (e.g., indexed by an explicit, fine parameter $\varepsilon$, periodic, etc.), then homogenization theory states the existence of a homogenized coefficient $A_0$ that replaces~$A$ in~\eqref{eq:model} and leads to a solution that presents the macroscopic behavior of the exact solution $u$, see~\cite{BlaLB23,Tar09}. In particular, a standard finite element method (as discussed above) may now be applied since problematic oscillations have been removed. There are constructive approaches to (approximately) obtain the coefficient $A_0$, such as, e.g., the heterogeneous multiscale method, see~\cite{EE03,EE05}. 

\noindent Other possibilities are computational strategies that evolve around modifying the space $V_h$ to include information about under-resolved scales. Corresponding strategies are, e.g., the multiscale finite element method~\cite{HouW97}, generalized multiscale finite element methods~\cite{EfeGH13,ChuEL18b}, multiscale spectral generalized finite element methods~\cite{BabL11,MaSD22,MaS22}, the localized orthogonal decomposition strategy~\cite{HenP13,MalP14,MalP20}, or gamblets~\cite{Owh17,OwhS19}. They often pose fewer assumptions on the structure of the coefficient~$A$. For a review article on computational strategies, see~\cite{AltHP21}. 

\noindent All these approaches have in common that they require additional pre-computations to obtain the homogenized coefficient or an appropriate discrete approximation space. These computations are the bottleneck of such strategies. Therefore, the benefit lies in scenarios where multiple right-hand sides $f$ in~\eqref{eq:model} are considered or in time-dependent cases in which $A$ does not change.

\subsection{Super-resolution}

Researchers recently paid attention to continuous image super-resolution (SR) for its practical and flexible image scaling for various displays. Local implicit image representation is one method that can map the coordinates and 2D features for continuous feature space interpolation. In this section, we will revisit these two topics in detail.

\noindent \textbf{Implicit neural representation.}
Implicit neural representation (INR) uses multi-layer perceptrons (MLP) to learn pixels or other signals from coordinates. It has been widely used in 3D shape modeling~\cite{shape_1,shape_2}, surface reconstruction~\cite{surface_1,surface_2}, novel view rendering~\cite{NERF,structure_1,structure_2}, etc. For instance, Mildenhall et al.~\cite{NERF} propose to map the camera poses and point clouds to the pixel values via a multi-layer MLP network. They use multiple-view images to optimize the network for implicit feature representation. Instead of using a voxel or point cloud, the implicit neural representation can 1) capture the fine details of scenes for photo-realistic reconstruction and 2) also emit complex 3D representation as a small number of differentiable network parameters.

\noindent INR can be regarded as an efficient, continuous, differentiable signal representation. Siren~\cite{siren} is one of the pioneering approaches that use periodic functions to build sinusoidal representation networks, such that it can solve complex boundary value problems. However, it suffers from signal noise and parameter variations. Therefore, a wavelet implicit neural representation (WIRE)~\cite{wire} is proposed to utilize Gabor wavelet activation functions to optimize time-frequency compactness. Its results on inverse problems of images and videos show that Gabor wavelets are a universally superior choice for nonlinearities. In the meantime, band-limited coordinate networks (BACON)~\cite{bacon} propose a similar solution as an analytical Fourier spectrum for multiscale signal representation. It initializes the frequencies of the sinusoidal layers of a multiplicative filter network (MFN)~\cite{mfn} within a limited bandwidth. Multiple layers of outputs are summed together to produce the overall output so that the network learns to fit the signal in a band-limited fashion. MIRE~\cite{mire} is proposed to use dictionary learning to match the activation of each layer for better signal representation. EVOS~\cite{evos} optimize the INR training via the modified evolution algorithm, which reduces $48\%\sim66\%$ training computation. ViiNeuS~\cite{viineus} proposes to use INR on a 3D surface representation. It learns the signed distance field from 2D images for LiDAR datasets.

\noindent In 3D processing, researchers propose to process 3D data by MLPs as unsigned distance~\cite{shape_1,shape_2}, signed distance~\cite{surface_1}, or occupancy field~\cite{surface_2}. For view rendering, NeRF~\cite{NERF} and Gaussian splitting~\cite{gaussians} also show great performance in real-time using implicit radiance field rendering. These advances in INR also inspire INR-based modeling of ordinary or partial differential equations. For example, Fourier neural operators~\cite{fno} parameterize the integral kernel directly in Fourier space with simple multiplications, which leads to a significant reduction of computational costs. Spherical Fourier neural operators \cite{sfno} generalize the Fourier transform, achieving stable autoregressive long-term forecasting.

\noindent \textbf{Continuous image super-resolution.}
Image and video super-resolution are ill-posed problems: One low-resolution (LR) image can lead to multiple plausible sharp and clean high-resolution (HR) images. The typical solution is to use paired LR-HR data to optimize restoration models via pixel-based loss functions. Since the seminal super-resolution work of Dong et al.~\cite{srcnn}, many image super-resolution approaches~\cite{sr_1,sr_2,sr3,srgan,swinir,stablesr,han,tdan,duf} have adopted an end-to-end neural network to learn the regression model for reconstruction. However, the existing state-of-the-art super-resolution methods can only handle fixed-scale upscaling. For out-of-distribution upscaling, they must be retrained or fine-tuned for optimal results. Continuous super-resolution, on the other hand, combines the recent developments of INR and super-resolution to come up with a grid-based latent space super-resolution. The local implicit image function (LIIF)~\cite{LIIF} is one of the pioneering strategies for continuous feature interpolation. That is, a grid-aware MLP network takes both pixel coordinates and global image features for feature expansion and then projects the super-resolved features to the actual pixel values. Such an operation can be further extended to video super-resolution. VideoINR~\cite{vinr} applies INR-based super-resolution to both space and time. Local texture estimators (LTE)~\cite{lte} improve them by estimating the dominant frequencies to explore the 2D Fourier space feature representation. Local implicit transformers (LITs)~\cite{clit} are used to integrate the attention mechanism and frequency encoding technique into the INR function. The cross-scale feature aggregation ensures a rich non-local representation. Given the recent success of diffusion models~\cite{stablesr} and attention to super-resolution~\cite{swinir}, combining INR and diffusion denoising models~\cite{diff_inr} shows superior performance. CiaoSR~\cite{ciaosr} proposes a multiscale attention-in-attention network for continuous super-resolution. Another direction is to use normalizing flow processes~\cite{normalize_flow}. The idea is to learn how to model the ill-posed problem as a super-resolution distribution. LINF~\cite{linf} learns a coordinate conditional normalizing flow model that can extract the local texture patches for global reconstructions. Most recently, NeurOp-diff~\cite{neuroop-diff} proposes to utilize a diffusion model based neural operator for remote sensing image super-resolution. HIIF~\cite{hiif} proposes a hierarchical encoding structure to replace MLP layers for continuous super-resolution. DDIR~\cite{ddir} combines the deformation field learning and implicit image function for real-world image super-resolution.

\subsection{Combining mathematical and super-resolution approaches}

Summarizing the above mathematical approaches, we recognize that numerical homogenization methods provably produce accurate solutions under specific model assumptions. However, these methods rely on a computational demanding offline phase for each fine-scale problem. That is, the fine-scale information correcting the coarse-scale solution needs to be calculated individually for each problem.

Machine learning-based super-resolution approaches learn the coarse to fine mapping for a vast class of fine-scale information. However, standard approaches do not respect the fine-scale structure of the individual problem; they rather learn a general fine to coarse mapping that is good for any fine-scale structure.

Our approach bridges this gap between mathematical \emph{numerical homogenization} and machine learning-based \emph{super-resolution}: It learns a super-resolution map that contains coarse-to-fine mappings that respect the individual fine-scale information. Thus, our model needs a simple and possibly inaccurate coarse-scale simulation and the fine-scale structure of the problem to generate a fine-scale approximation to the solution.
As this approach is very general and efficient, it cannot provide results of the quality of other specifically constructed numerical homogenization tools that invest additional resources for every coefficient separately. Nonetheless, we achieve very good results in a macroscopic sense, which is already very important for multiscale problems. In particular, the results are much better than straight-forward coarse simulations.
However, there is (as almost always for machine learning models) no mathematical proof that the result satisfies mathematical accuracy bounds. While one could check the quality of the obtained solution a posteriori, which is a future research prospect, we rather want to provide an intuitive justification for our approach:

The model~\eqref{eq:model} can be viewed as a description of (slow) fluid flow through soil. 
The coarse-scale approximation discards fine properties of the flow and is therefore inaccurate. The missing piece is still encoded in the actual soil structure.
Since our approach uses both pieces of information (coarse flow and fine-scale structure of porous matrix) and learns to transfer simulation results from a coarse to a fine scale, we can hope to obtain good solutions. This is experimentally justified in Section 4.

Building on this motivation, we further highlight how our method differs from existing super-resolution approaches. Unlike prior works that are limited to fixed-scale super-resolution, our proposed NH-CSR achieves scale-arbitrary super-resolution for both in-distribution and out-of-distribution cases within a single unified model. At its core, NH-CSR employs a multiscale implicit image function, which surpasses existing continuous super-resolution methods in interpolating missing sub-pixels across varying scales. Furthermore, unlike traditional implicit neural representations that require per-signal optimization, NH-CSR adopts a meta-learned, coordinate-based architecture that enables general image representation with a single shared model.

\begin{figure*}[t]
	\begin{center}
		\centerline{\includegraphics[width=\textwidth]{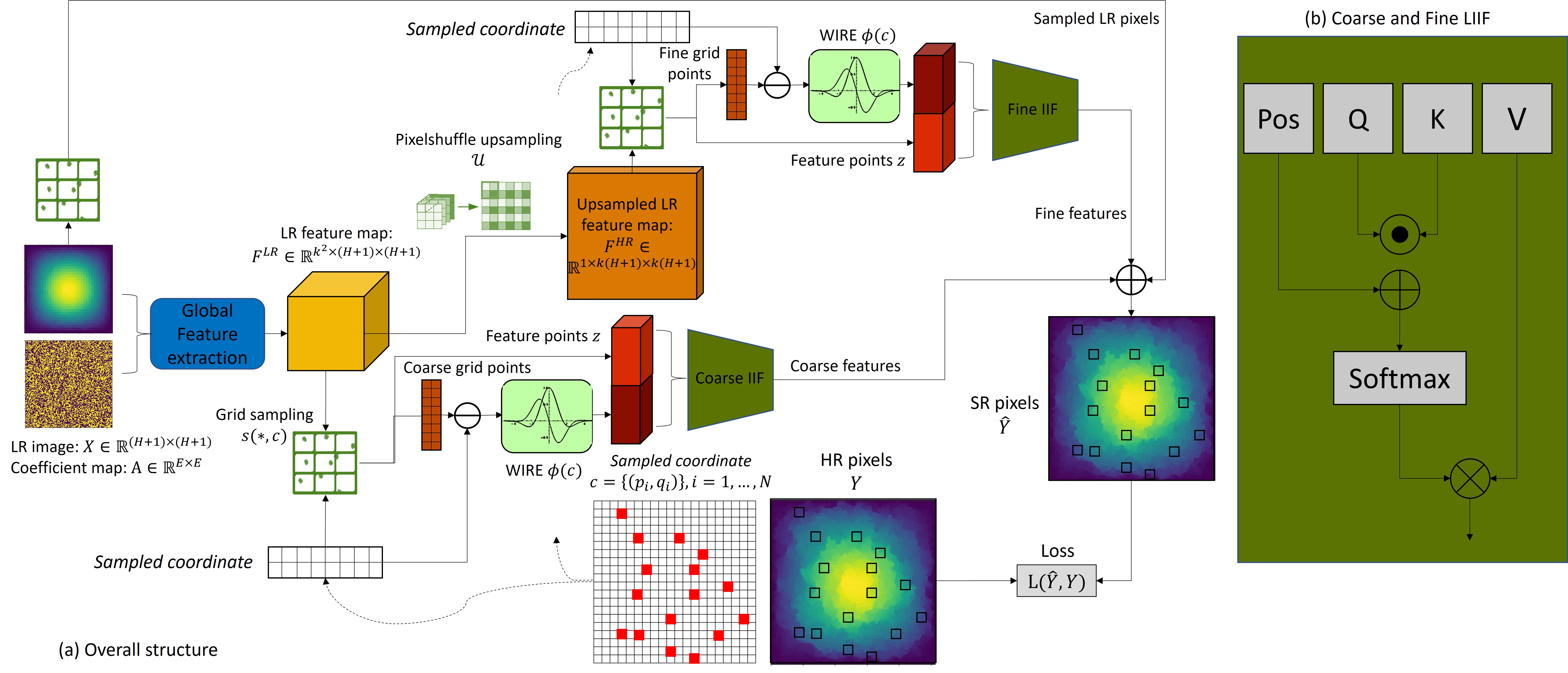}}
        \caption{\small \textbf{The proposed super-resolution model for numerical homogenization.} (a) shows the overall architecture and (b) shows the detailed operator of LIIF used in coarse and fine LIIF modules. Given a coarse fine element solution represented by $X$ and the corresponding coefficient $A$, we learn a conditional feature extractor using off-the-self neural networks. The extracted feature is then sampled based on the random coordinate map $c$, so that we obtain sampled feature points and grid points. Coarse and fine implicit image functions (IIF) take as input the paired feature and grid points learned by the proposed WIRE operation to form the regression model for super-resolution. The results of coarse- and fine-grained pixel regression are summed together with the sampled LR pixels to obtain the final super-resolved pixels.}
		\label{fig:network}
	\end{center}
\end{figure*}

\section{Approach: numerical homogenization by super-resolution}
\label{Approach}

The aim of our approach is to only compute a (not necessarily reasonable) coarse approximation based on the finite element method that is then upscaled to a finer scale by a continuous super-resolution approach without using sophisticated mathematical numerical homogenization tools. From knowledge of a coarse finite element solution (which is cheap to obtain but in general not even a good approximation in a macroscopic sense) and the coefficient $A$, we aim to recover a reasonable solution on a finer scale which contains microscopic information as well. 
To connect this setup to image-based super-resolution strategies, we represent the coarse finite element solution as a low-resolution (LR) pixel image, the fine-scale solution by a high-resolution (HR) image, and the coefficient as an image on an intermediate scale. 
We aim to upscale solutions independently of the scale difference: starting from a certain coarse approximation on the scale $h$ with some coefficient that is piecewise constant on a grid consisting of boxes with side length $\varepsilon$, our method should be able to recover an upscaled solution on a chosen scale $\delta$, where typically $h > \varepsilon > \delta$. 
This approach may be understood as a numerical homogenization strategy in the spirit of the above-mentioned multiscale approaches, but it acts a posteriori. Its idea of upscaling is more abstract than that of traditional numerical homogenization approaches and does not boil down to modifying basis functions individually. In particular, our approach is not necessarily local and different parts of the image and also features potentially learned from other images are taken into account for a global prediction. 
This is achieved by a continuous super-resolution model, which is trained on a set of coefficients and corresponding coarse and fine solutions. More details are given in the following subsection.  

\noindent Our numerical homogenization strategy via continuous super-resolution (\net) is illustrated in Figure~\ref{fig:network}. We denote by $X\in\mathbb{R}^{(H + 1)\times (H + 1)}$ the input low-resolution finite element data, where {$H + 1$} is the dimensions of its height and width {and $H = 1/h$}. The corresponding coefficient $A\in\mathbb{R}^{E\times E}$ {with $E = 1/\varepsilon$} is also taken as an additional input for the super-resolution model. Note that $A$ was previously defined as a function. However, since it is assumed to be piecewise constant, it may be represented as a pixel image as well. The overall neural network is a coefficient-guided continuous super-resolution model $f(X, A, \alpha)$ to produce an $\alpha$-times super-resolved finite element solution $\hat{Y}\in\mathbb{R}^{{(\alpha H + 1)}\times {(\alpha H + 1)}}$, where $\alpha \in \mathbb{R}^{+}$. As shown in Figure~\ref{fig:network}(a), the complete network consists of three main components: 1) global feature extraction (blue block), 2) learnable Gabor wavelet-based grid sampling (The WIRE in light-green blocks), and 3) Multiscale implicit image function (MS-IIF) (olive-green blocks, Figure~\ref{fig:network}(b)).

\subsection{Global feature extraction.} Intuitively, to achieve arbitrary-scale super-resolution, we first need to have a continuous data representation. The coarse finite element data and the coefficient map are discrete with different data ranges. Hence, we propose a conditional feature extraction using an off-the-shelf image processing backbone to extract the dense feature representation as $F=g(d(X), A)$ (blue block in Figure~\ref{fig:network}(a)). We first upsample the coarse finite element data via the bicubic interpolation $d(\cdot)$ to the same size as the coefficient map~$A$
and then stack them together as input to learn the dense features $F$ via several convolution operations $g$. 
Note that we use the off-the-shelf networks as backbones without using their pre-trained parameters as the finite element data are generated randomly without strong semantic meaning as for natural photos. Using more complex networks can gain some improvements. We will discuss these in Section~\ref{Experiments}.

\subsection{Gabor wavelet based implicit image function.} For continuous super-resolution, we are interested in learning a local implicit image function (LIIF)~\cite{LIIF},

\begin{small}
\begin{equation}
{J} = f_\theta(F, r),
\label{eq:liif}
\end{equation}
\end{small}

\noindent where $J$ is the predicted data and $F$ is the global feature map obtained from the previous step. Further, $r$ is the 2D coordinate in the continuous 2D domain and $f_\theta$ is an $M$-layer MLP network. We assume that $f_\theta$ can be shared by all the images, hence it learns the local feature representation by assigning local feature points to their 2D coordinates,

\begin{small}
\begin{equation}
{J} = f_\theta(s(F, c), s(r, c)-c),
\label{eq:grid}
\end{equation}
\end{small}

\noindent where $s$ is the bilinear grid sampling operation, $r$ represents the 2D coordinates of the global feature map, and $c$ are the target 2D coordinates. We first apply spatial sampling to the global feature map and its coordinates to localize the neighborhood features, then we use MLP to reconstruct the missing data. Here, we use the coordinate residues $s(r, c)-c$ to emphasize the distances between nearest neighbors for prediction. 

The term in~\eqref{eq:grid} describes a grid-based regression model. Using MLP with ReLUs is biased towards learning low-frequency content~\cite{LIIF}. Inspired by the recent development of positional encoding~\cite{attention} and Fourier feature mapping~\cite{NERF}, implicit neural representation shows much better results on high-frequency information prediction. Consider an implicit neural representation function $w_\theta:\mathbb{R}^{D_i}\rightarrow\mathbb{R}^{D_o}$ mapping input coordinates to $D_o$ output dimensions with learnable parameters $\theta$. The goal is to map the coordinates to another domain of interest. In general, based on a coordinate vector $y_{m-1}$ in the $(m-1)$-th layer, the consecutive coordinate vector is computed by

\begin{small}
\begin{equation}
y_m = w_\theta(y_{m-1}) = \sigma(W_m y_{m-1} + b_m),
\label{eq:gabor}
\end{equation}
\end{small}

\noindent where $\sigma$ is the nonlinear activation function and $W_m$, $b_m$ are the weights and biases for the $m$-th layer. In previous works~\cite{lte,clit}, learnable $M$-layer MLPs are utilized for coordinate expansion, including the periodic function $\sigma(x)=\sin(\omega x)$ and the plain ReLU function $\sigma(x)=\max(\omega x, 0)$, where~$\omega$ denotes the learning parameter for the MLP layer. The problem is that they do not offer an explicit multi-resolution representation. In our case, the finite element data naturally preserves a multi-resolution ability, and therefore we propose continuous Gabor wavelets $\phi$ for the coordinate mapping (light-green blocks in Figure~\ref{fig:network}(a)), $\phi(x; \omega_0, s_0)=\cos(\omega_0 x)e^{-(s_0x)^2}$. It is interesting to point out that when $s_0=0$, we have $\phi(x; \omega_0, s_0)=\cos(\omega_0 x)$ which is similar to SIREN in~\cite{siren}. When $\omega_0=0$, we have $\phi(x; \omega_0, s_0)=e^{-(s_0x)^2}$, which is the Gaussian nonlinearity in~\cite{gaussian_cor}. We apply two layers of such Gabor wavelet nonlinearities to the 2D coordinates, such that we can rewrite~\eqref{eq:grid} as

\begin{small}
\begin{equation}
J = f_{\theta}(s(F, c), \phi\left(s(r, c)-c\right)),
\label{eq:gabor_grid}
\end{equation}
\end{small}

\noindent where $J$ is the output feature vector. Furthermore, we modify the attention mechanism to expand the local feature aggregation. Specifically, we reformalize~\eqref{eq:gabor_grid} as a non-local filtering process (Figure~\ref{fig:network}(b)),

\begin{small}
\begin{equation}
\begin{matrix}
\!\begin{aligned}
& Q = s(W_\mathrm{Q}(z)),\  K = s(W_\mathrm{K}(z)), \ V = s(W_\mathrm{V}(z)), \\
& J=f_\theta\left[\mathrm{softmax}\left(\frac{QK^T}{\sqrt{D}}+\phi(s(r, c)-c)\right)\times V\right],
\label{eq:attn}
\end{aligned}
\end{matrix} 
\end{equation}
\end{small}

\noindent where $Q, K, V$ are the query, key, and value to compute the attention. Further, $D$ is the channel dimension of the local feature embedding. It first calculates the inner product of $Q$ and $K$, adds the relative positional bias obtained by a Gabor wavelet based coordinate encoding to the result, and produces an attention matrix normalized by the softmax operation. Intuitively, \eqref{eq:attn} is the weighted version of~\eqref{eq:gabor_grid}, which learns the feature correlations and aggregates similar patterns for efficient image reconstruction. 

\subsection{Multiscale implicit image function (MS-IIF).} A key challenge in continuous super-resolution is how to faithfully interpolate across scales: low-frequency components must be preserved to maintain global consistency, while high-frequency details need to be accurately reconstructed for sharp visualization. A single-scale representation often struggles to balance these two aspects, leading either to overly smooth results (if biased toward low frequencies) or to noisy artifacts (if biased toward high frequencies). To address this, we propose a grid-based \textit{multiscale implicit image function} (MS-IIF) (coarse and fine LIIF modules in Figure~\ref{fig:network}(a)), which explicitly disentangles coarse- and fine-scale interpolation. It takes the extracted feature maps $J\in\mathbb{R}^{k^2\times {(H+1)\times (H+1)}}$ (where $k$ is the number of channels) and expands the values to a higher dimension. In the feature domain, we want to upscale the features to the desirable size of $k(H + 1)\times k(H + 1)$. To achieve that, we first generate a coarse grid $r^\mathrm{LR}\in\mathbb{R}^{2\times {(H+1)\times (H+1)}}$ and a fine grid $r^\mathrm{HR}\in\mathbb{R}^{2\times {k(H+1)\times k(H+1)}}$. Then we apply a coarse and a fine data regression, that is,

\begin{small}
\begin{equation}
\begin{matrix}
\!\begin{aligned}
& F^\mathrm{LR}=f_{\theta}^\mathrm{LR} [s(J, c), \phi(s(r^\mathrm{LR}, c)-c)], \\
& F^\mathrm{HR}=f_{\theta}^\mathrm{HR} [s(\mathcal{U}(J), c), \phi(s(r^\mathrm{HR}, c)-c)], \\
& F=\mathrm{MLP}(F^\mathrm{LR} + F^\mathrm{HR}).
\label{eq:multi}
\end{aligned}
\end{matrix} 
\end{equation}
\end{small}

\noindent Using the pixelshuffling operator $\mathcal{U}$ to upsample the original global feature maps $J$, we obtain the reconstructed feature map~$F$ as the sum of coarse-scale ($F^\mathrm{LR}$) and fine-scale ($F^\mathrm{HR}$) grid-sampled features. Finally, we use one more MLP layer $\mathrm{MLP}$ to obtain the final reconstructed finite element data points.

\subsection{Losses $\mathrm{L}(Y, \hat{Y})$}. To train our model, we combine an $L1$ loss and our proposed stochastic cosine similarity for optimization. The $L1$ loss is given by $\mathrm{L}_1(Y, \hat{Y}) = \frac{1}{B}\sum_{i=1}^B \|{Y_i}-\hat{Y}_i\|^1_1$, where $B$ is the batch size. Note that here $Y_i$ and $\hat{Y}_i$ are the sampled pixels from predicted and ground truth data, respectively. The problem with this loss is that it takes each data point individually to measure the differences. For a 2D finite element solution, the local region shows similar continuous patterns that need to be collectively optimized. In order to supervise the network to learn data points as a whole set, we get inspiration from S3IM~\cite{s3im} to propose a stochastic cosine similarity score (SCS) $\mathrm{L}_{\mathrm{scs}}$ as follows,

\begin{small}
\begin{equation}
\mathrm{L}_\mathrm{scs}(Y, \hat{Y}) = \frac{1}{B(H+1)} \sum_{b=1}^B \sum_{t=1}^{H+1} \frac{\mathcal{P}^t_b(\hat{Y}) \times \mathcal{P}^t_b(Y)}{\|\mathcal{P}^t_b(\hat{Y})\| \|\mathcal{P}^t_b(Y)\|}.
\label{eq:scs}
\end{equation}
\end{small}

\noindent The SCS can be summarized in three steps:

\begin{enumerate}
\item We randomly sample a batch size of $B$ input data points, each datum is of dimension $1\times {(H+1) \times (H+1)}$. The corresponding ground truth and model output $Y$ and $\hat{Y}$  are of dimension $B\times {(H+1)\times H+1}$, where the ground truth data points are randomly sampled from the original images. We randomly sample $H+1$ points from $Y$ and $\hat{Y}.$ 

\item We compute the cosine similarities  $\mathcal{P}^t_b(Y)$ and $\mathcal{P}^t_b(\hat{Y})$.

\item We repeat steps (1) and (2) $H+1$ times and average the $H+1$ estimated stochastic cosine similarities to obtain the final SCS.
\end{enumerate}

\noindent We avoid using edge-guided loss functions, e.g., SSIM, because they focus more on the edges and may maximize the high-frequency details. For finite element super-resolution, we prefer data fidelity over visual alignment. To train the proposed \net, we have the total loss as $\mathrm{Loss} = \mathrm{L}_1(\hat{Y}, Y)+\lambda \mathrm{L}_\mathrm{scs}(\hat{Y}, Y)$. The weighting parameter~$\lambda$ balances the two loss terms. 

\subsection{Physical re-interpretation}

From a modeling perspective, our approach tries to achieve the impossible: it generates/guesses the fine-scale behavior (oscillations) of a solution from a coarse (smoothened) solutions. From an abstract point of view, this can be interpreted as the correct generation of turbulence or any fine-scale behavior typically ignored (averaged) in numerical simulations. To this end, we solely exploit tools from image super-resolution without using the actual physics of the problem (these physics are properly learned by our approach without building them into the model).

\section{Experiments}
\label{Experiments}
\subsection{Implementation Details}
\begin{figure}[t]
	\begin{center}
		\centerline{\includegraphics[width=\columnwidth]{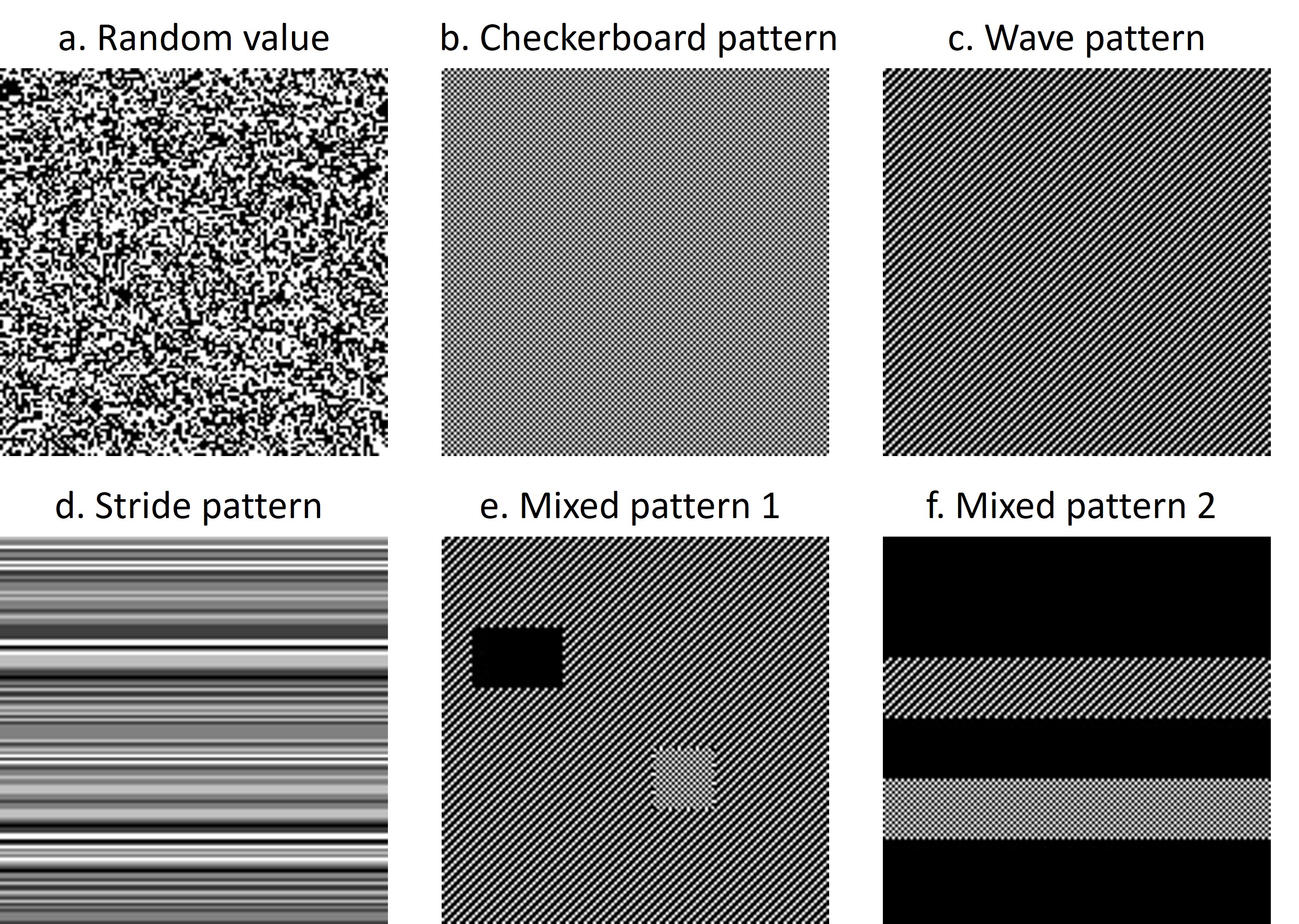}}
        \caption{\small \textbf{Examples of coefficient maps used for finite element data generation.} We show six examples of different coefficients. The visualization is done by normalizing the coefficient value to [0, 255].}
		\label{fig:coef}
	\end{center}
\end{figure}

\noindent \textbf{Datasets.}
To generate sufficient data samples for model training and analysis, we randomly generate coefficient maps (independently and identically distributed random values in~$\{1,100\}$ on a grid of size~$128\times 128$ (which corresponds to $\varepsilon = 1/128$)) and use a standard finite element code to obtain the corresponding low-resolution and high-resolution finite element solutions as data (we interpret nodal values as pixels). The resolution of LR data is $33\times 33$ (which corresponds to $h = 1/32$), the HR data is $(32\alpha + 1)\times (32\alpha + 1)$, where $\alpha=2,4,8,16,32,64$. The scaling of alpha by multiples of two and the respective size of the images is classical for finite elements, where the mesh is successively refined by halving the mesh size (i.e., the side length of pixels). For the training, we randomly generate binary coefficient maps as mentioned above. For testing, we have two scenarios: 


\begin{enumerate}
\item \textbf{Random testing} generates another set of paired LR-HR data using random binary coefficient maps.
\item \textbf{Deterministic testing} generates new paired LR-HR data based on deterministic coefficient maps. We focus two specific patterns (a wave and a checkerboard structure), and we mix them to form 13 new coefficient maps, see Figure \ref{fig:coef}.
\end{enumerate}

\noindent We show six examples of different coefficient patterns in Figure~\ref{fig:coef}: (a) is an example of a random coefficient map that we use to generate paired training data, (b) to (d) are three deterministic coefficients generated from checkerboard, wave, and striding patterns, and (e) and (f) are two examples of mixed coefficient maps where we randomly mix three deterministic coefficients for testing. 

\noindent We consider two super-resolution scenarios:
\begin{enumerate}
\item \textbf{In-distribution SR:} We train our model for $16$-times super-resolution and test the model under the same super-resolution factor. Since the model is optimized using paired LR-HR data under $16$-times super-resolution, the model learns the fixed mapping correlations for super-resolution. 
\item \textbf{Out-of-distribution SR:} We train our model for $16$-times super-resolution and test the model under other super-resolution factors, including $\alpha=2,4,8,32,64$. The corresponding output data size is $65 \times 65, 129 \times 129, 257 \times 257, 1025 \times 1025, 2049 \times 2049$, respectively.
\end{enumerate}

\noindent \textbf{Parameter setting.}
The global feature extractor can be any off-the-shelf image encoder. We use residual dense networks (RDN)~\cite{rdn} and Swin~\cite{swinir}. In the following discussion, we use RDN as the default encoder as it balances quality and efficiency. We train \net\ using the Adam optimizer with the learning rate $10^{-4}$. The learning rate is halved after 50.000 iterations. The batch size is set to 32, and \net\ is trained for 100.000 iterations (about 4 hours) on a PC with one NVIDIA V100 GPU using the PyTorch deep learning platform. 

\begin{table*}[b]\centering
\caption{\textbf{Comparison with state-of-the-art methods.} We report the results on the test datasets of the arbitrary upscaling scenarios. PSNR and MSE ($10^{-3}$) are used for comparison.} 
\renewcommand\arraystretch{1.4}
\resizebox{\textwidth}{!}{
\begin{tabular}{c|cccccccc|cccc}
\toprule
\multirow{3}{*}{Method} & \multicolumn{8}{c|}{In-distribution} & \multicolumn{4}{c}{Out-of-distribution} \\ \cline{2-13} 
 & \multicolumn{2}{c}{2x} & \multicolumn{2}{c}{4x} & \multicolumn{2}{c}{8x} & \multicolumn{2}{c|}{16x} & \multicolumn{2}{c}{32x} & \multicolumn{2}{c}{64x} \\
 & MSE$\downarrow$ & SSIM$\uparrow$ & MSE$\downarrow$ & SSIM$\uparrow$ & MSE$\downarrow$ & SSIM$\uparrow$ & MSE$\downarrow$ & SSIM$\uparrow$ & MSE$\downarrow$ & SSIM$\uparrow$ & MSE$\downarrow$ & SSIM$\uparrow$ \\ \midrule
MetaSR~\cite{METASR} & 0.5717 & 0.7555 & 0.8902 & 0.8571 & 1.0274 & 0.9004 & 1.3115 & 0.9048 & 1.6404 & 0.9157 & 1.9924 & 0.9201 \\
LIIF~\cite{LIIF} & 0.5912 & 0.7493 & 0.8742 & 0.8574 & 0.9523 & 0.9222 & 1.2045 & 0.9329 & 1.5122 & 0.9291 & 1.8445 & 0.9198 \\
LTE~\cite{lte} & 0.5774 & 0.7070 & 0.9044 & 0.8439 & 0.9554 & 0.9199 & 1.1640 & 0.9329 & 1.4442 & 0.9005 & 1.7445 & 0.8948 \\
IDM~\cite{idm} & 0.6152 & 0.7110 & 0.9125 & 0.8442 & 1.0027 & 0.9145 & 1.2204 & 0.9223 & 1.5548 & 0.9012 & 1.8667 & 0.8910 \\
LIT~\cite{clit} & 0.5719 & 0.7524 & 0.8654 & 0.8578 & 0.9246 & 0.9234 & 1.1576 & 0.9335 & 1.4972 & 0.9300 & 1.8056 & 0.9217 \\
\cellcolor{mistyrose}{Ours} & \cellcolor{mistyrose}{0.5706} & \cellcolor{mistyrose}{0.7586} & \cellcolor{mistyrose}{0.8340} & \cellcolor{mistyrose}{0.8627} & \cellcolor{mistyrose}{0.8848} & \cellcolor{mistyrose}{0.9243} & \cellcolor{mistyrose}{1.0743} & \cellcolor{mistyrose}{0.9346} & \cellcolor{mistyrose}{1.3825} & \cellcolor{mistyrose}{0.9310} & \cellcolor{mistyrose}{1.7020} & \cellcolor{mistyrose}{0.9227} \\ \bottomrule
\end{tabular}
}
\label{tab:sota}
\end{table*}

\noindent \textbf{Metrics and evaluation.}
Mean squared errors (MSE) and mean absolute errors (MAE) measure the average pixel differences between ground truth and estimation. {Structural similarities} (SSIM)~\cite{ssim} measures the structural closeness of ground truth and estimation. To measure the model complexity, we use FLOPs, memory, and the number of parameters to compare the computational efficiency.

\noindent \textbf{Experimental design and robustness.}
Our experimental setup is grounded in the theory of numerical homogenization, a general mathematical framework that exploits fine-scale solution features (often referred to as oscillations) within coarse-scale simulations. In this context, testing on finite element data with randomly generated coefficients is a standard and theoretically justified practice, as it effectively captures a broad spectrum of possible variations in the underlying medium. These randomized settings are sufficient to evaluate both the robustness and the generalization capability of the proposed method, thereby addressing the main concerns typically raised about experimental validity. The influence of real-world noisy data is not addressed within this manuscript, as its influence is usually analyzed by running many (similar) numerical experiments, which all could be enhanced using our approach.

\begin{figure*}[ht]
	\begin{center}
		\centerline{\includegraphics[width=\textwidth]{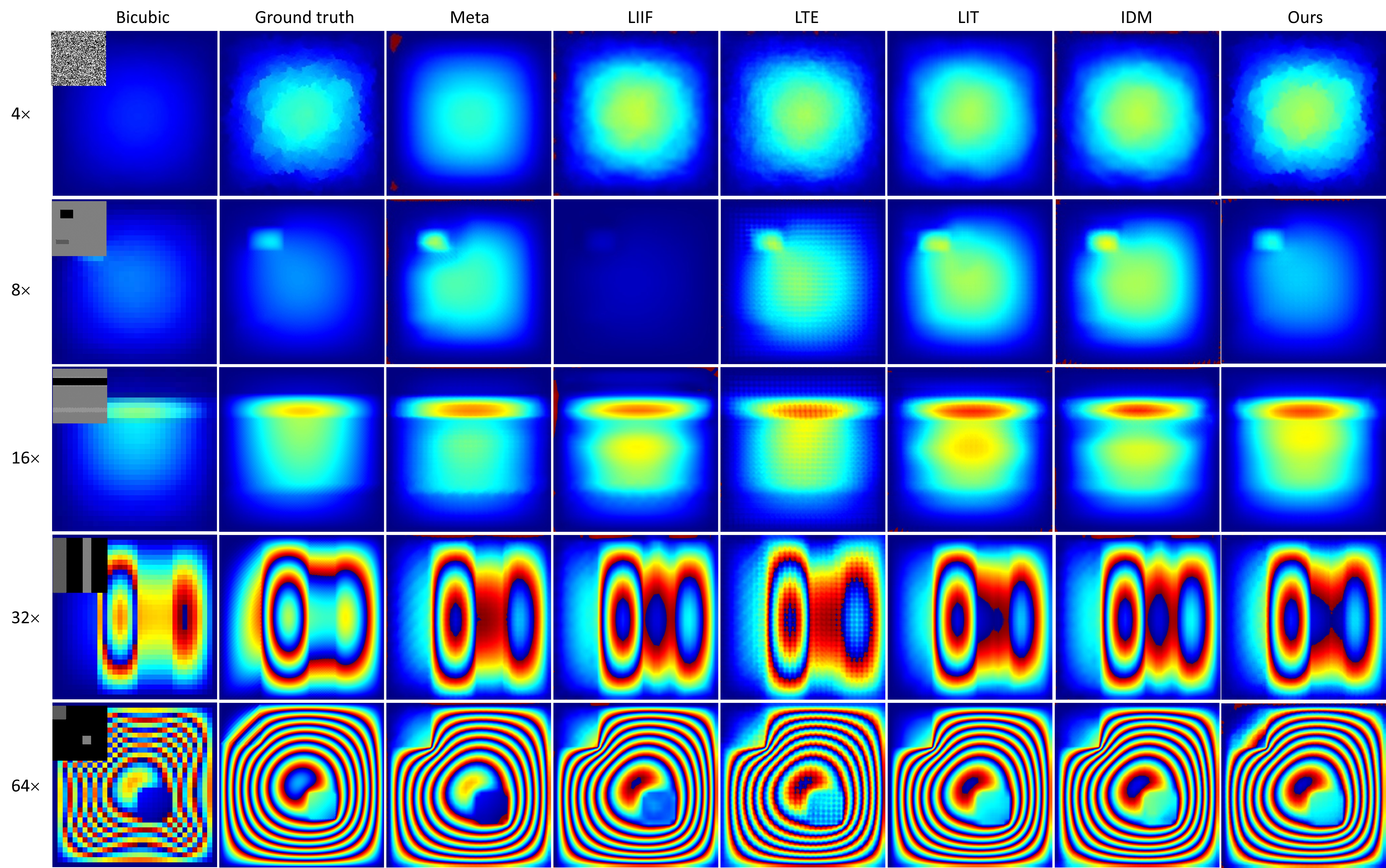}}
        \caption{\small \textbf{Visual comparison of the coarse finite element solutions.} We use the proposed method to apply multiscale upsampling to the coarse finite solution, both in-distribution (upsampling factor no larger than 16) and out-of-distribution (upsampling factor larger than 16) scenarios. Ours achieves the best visualization among others. As before, the pixel values are multiplied by a factor 8 and shifted to a value range between 0 and 255 to better visualize differences. More precisely, the multiplied values are taken modulo 255.}
		\label{fig:sota_mix}
	\end{center}
\end{figure*}

\subsection{Comparison of the continuous super-resolution with state-of-the-art methods}\label{sota}
To show the efficiency of our proposed method, we compare it with four state-of-the-art continuous super-resolution methods: Meta-SR~\cite{METASR}, LIIF~\cite{LIIF}, LTE~\cite{lte}, LIT~\cite{clit} and IDM~\cite{idm}. All approaches are reimplemented using the publicly available codes. We retrained their model on our finite element dataset for fair comparison. In Table~\ref{tab:sota}, we report the results of using different methods on both in-distribution and out-of-distribution super-resolution scenarios. We can see that ours achieves the best performance in both MSE and SSIM across different upsampling factors. Note that we use MSE to measure the data fidelity between estimation and the ground truth. SSIM is used to measure the structural closeness which indicates whether super-resolution can preserve the semantic patterns. We can also observe that LIIF and LIT are the most competitive approaches. Meta and LTE have lower performance than others, and IDM also does not perform well as it is built on the pre-trained autoencoder~\cite{idm}. Compared with others, our model shows consistent significant improvements on MSE by $0.01 - 0.3$ and SSIM by $0.01 - 0.03$. 

\noindent For better comparison, we project super-resolved finite element data to the RGB domain for visualization. To do so, we first normalize the 2D finite element data to [0, 255] and then multiply the value by 8 in order to highlight the subtle numerical changes. Note that values above 255 are shifted to a value range between 0 and 255 by taking the values modulo 255. This leads to a highlighting of local differences. Finally, we use the \textit{JET} colormap to expand the data to the RGB space. The visual comparison is shown in Figure~\ref{fig:sota_mix}. Here we show the results of the super-resolution strategy based on coarse finite elements generated by different coefficients (e.g., (e) and (f) in Figure~\ref{fig:coef}). In the first row, we have a setting generated by a random coefficient, where we observe that our model shows clearer curvatures compared to other models. In the second row, LIIF fails to predict the patterns, and except ours, other models show an exaggerated behavior in the center regions. The third row shows that our model predicts best the general geometric configuration of the HR solution. Finally, in rows 4 and 5 our method better predicts the blue channel between the two circles in row 4 and the upper left edge in row 5. 

\begin{figure*}[ht]
	\begin{center}
		\centerline{\includegraphics[width=\textwidth]{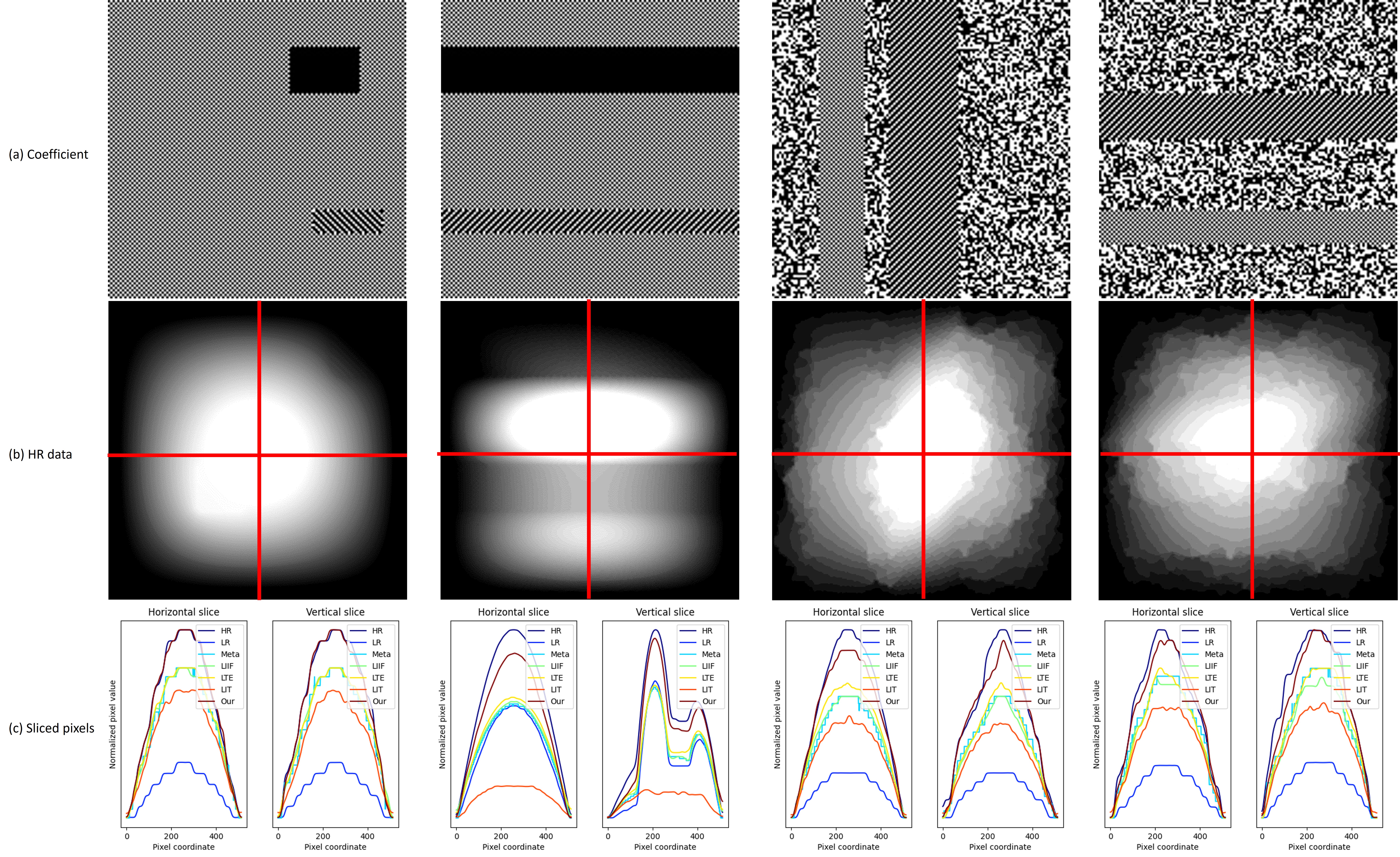}}
        \caption{\small \textbf{Plots of sliced pixels from different super-resolution approaches.} To better visualize the individual pixel differences, we use different methods to apply 16-times super-resolution on one coarse finite element data, then we slice the middle row and column of pixels (indicated by the red lines) and plot their values by coordinates. We can see that our approach better aligns with the HR data.}
		\label{fig:slice_pixel}
	\end{center}
\end{figure*}

\noindent Note that the visualization of actual finite element data distorts the lower and upper bounds of the data as mentioned above. Therefore, a good way to further visualize the super-resolution is to plot certain horizontal and vertical slices and compare the pixel values along these slices. We can see from Figure~\ref{fig:slice_pixel} that our method (dark brown) better aligns with the ground truth data (dark blue), especially in the first example, while others show different degrees of misalignment. On the other hand, the other competitive methods (LTE in yellow color, LIIF in green color, and LIT in red color) struggle in these four cases to get close to the ground truth. Our method works extremely well in capturing results for random coefficients (that may only attain two different values per pixel) and even beyond (cf.~particularly the first example with a deterministic pattern). Note that in other examples, our method does not perform as well as for the first and the random cases. This is related to the fact that the coefficient maps that underlie these cases include patterns that are included in the training phase, where only random patterns were considered. Additionally, those patterns that generate problems contain \emph{channels} through the whole domain yielding a different physical behavior from what the neural network has seen and also requiring highly specified tools in classical numerical homogenization. 
Nevertheless, our results are still much better than those of all the other super-resolution approaches, which gives great hope in achieving even better results if more involved training sets are considered beyond the random case.

\noindent Finally, we are also interested in the frequency analysis on the super-resolution results to see the frequency interpolation. The natural way to describe the high-frequency behavior is through a Fourier transform. To do so, we apply a radial power spectral density (RAPSD), which is based on the assumption that by slicing at any angle of the frequency spectrum, we are able to identify changes in the magnitude spectra. If the model is able to predict the missing high-frequency components, the corresponding energy should match. As shown in the left column of Figure~\ref{fig:rapsd},  we apply RAPSD to out-of-distribution (32-times and 64-times) SR results in the following steps: 

\begin{figure}[H]
	\begin{center}
		\centerline{\includegraphics[width=\columnwidth]{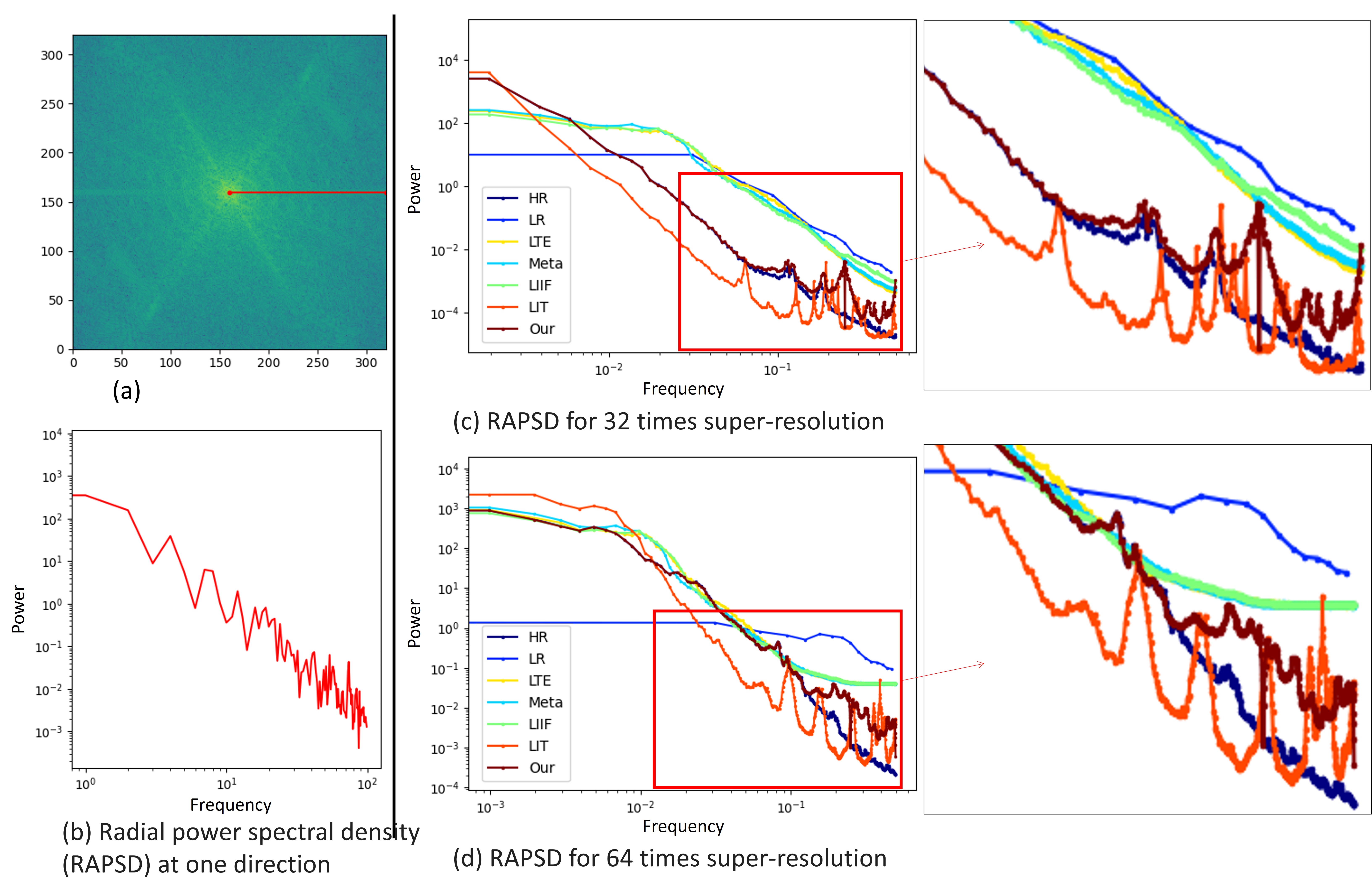}}
        \caption{\small \textbf{Frequency analysis of different super-resolution approaches.} We use the Fourier transform to convert super-resolved data to the frequency domain, and then we use radial power spectral density (RAPSD) to average the spectrum along different angles. Finally, we use the log-log plot to show the correlations between power and frequency. We can see that our approach can better preserve the power across different frequency bands.}
		\label{fig:rapsd}
	\end{center}
\end{figure}

First, we convert the SR data to the Fourier domain and slice the frequency spectrum in the clockwise direction to sample the spectral magnitude changing along frequencies ((b) in Figure~\ref{fig:rapsd}). Then we average the power-vs-frequency results on all testing data and plot it in~(c) and~(d) in Figure~\ref{fig:rapsd}. We can see that ours better aligns the trend of energy changes along the frequency bands. On the other hand, LIT introduces many spikes along the high-frequency bands. LIIF, Meta and LTE (green, blue and yellow curves) align first on the low-frequency domain and then introduce errors in the middle- and high-frequency domains. 

\begin{figure*}[ht]
	\begin{center}
		\centerline{\includegraphics[width=\textwidth]{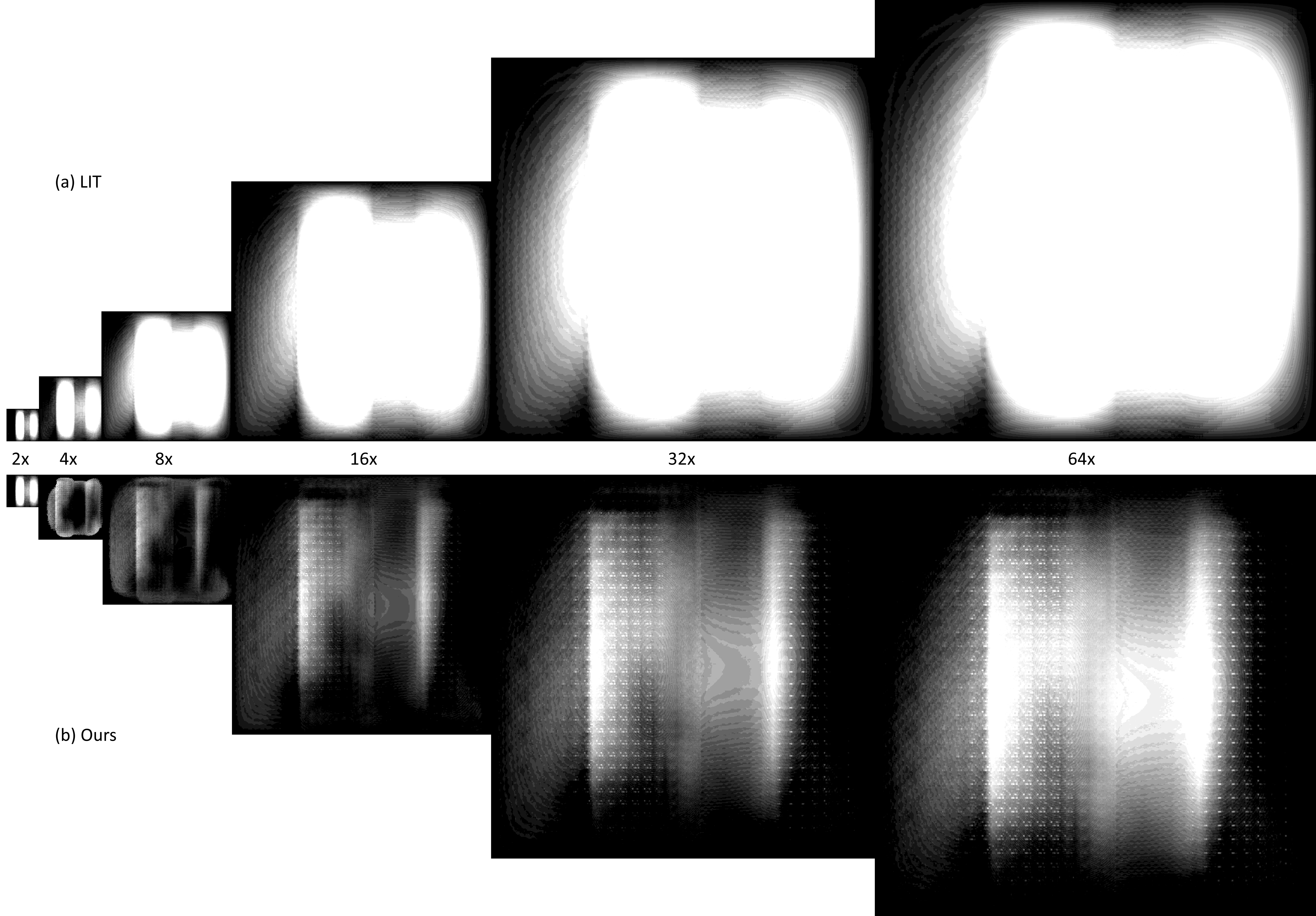}}
        \caption{\small \textbf{Visualization of continuous super-resolution to finite elements.} We show the residual image between LIT and ours on various upsampling factors. For better visualization, we multiply the residual values by 20 times.}
		\label{fig:csr}
	\end{center}
\end{figure*}

\subsection{Ablation studies and analysis}
For a further analysis of our approach, we conduct several ablation studies to demonstrate the efficiency of key modules.

\noindent \textbf{Continuous super-resolution.} We first need to analyze whether the proposed model can consistently produce super-resolved finite element approximations. From Table~\ref{tab:sota}, we can find that the most competitive approach for us is LIT. We therefore compare our results to LIT by visualizing the residual map between SR and HR data. In Figure~\ref{fig:csr}, we show the residual maps for different upsampling factors. Whiter means a larger residual in that context. For all upsampling factors, we perform better than LIT.

\noindent \textbf{Gabor based grid sampling.} We are also interested in the effect of using Gabor-based grid sampling for the coordinate encoding. To validate its effect, we compare it with different encoding schemes, including plain encoding without learnable parameters, sinusoid~\cite{siren}, IPE~\cite{NERF}, and LTE~\cite{lte}. In Table~\ref{tab:gabor}, we can see that our method achieves the best performance in terms of MSE, SSIM, and MAE.

\begin{table}[H]\centering
\caption{\textbf{Ablation study on the key modules.} We report the results on the test datasets of 4-times and 8-times upsampling, including MSE ($10^{-3}$), MAE ($10^{-3}$) and SSIM.} 
\renewcommand\arraystretch{1.4}
\resizebox{\columnwidth}{!}{
\begin{tabular}{cc|ccccc}
\toprule
\multicolumn{2}{c|}{Positional encoding} & Plain & Sinusoid & IPE & LTE & \cellcolor{mistyrose}{Gabor} \\ \midrule
\multirow{3}{*}{4x} & MSE$\downarrow$ & 1.002 & 0.989 & 0.926 & 0.974 & \cellcolor{mistyrose}{0.834} \\
 & SSIM$\uparrow$ & 0.851 & 0.857 & 0.855 & 0.857 & \cellcolor{mistyrose}{0.863} \\
 & MAE & 21.76 & 21.44 & 21.05 & 21.17 & \cellcolor{mistyrose}{20.09} \\ \midrule
\multirow{3}{*}{8x} & MSE$\downarrow$ & 1.124 & 1.110 & 1.057 & 0.985 & \cellcolor{mistyrose}{0.885} \\
 & SSIM$\uparrow$ & 0.917 & 0.917 & 0.918 & 0.919 & \cellcolor{mistyrose}{0.924} \\
 & MAE$\downarrow$ & 17.84 & 17.46 & 17.37 & 17.28 & \cellcolor{mistyrose}{16.93} \\
 \bottomrule
\end{tabular}%
}
\label{tab:gabor}
\end{table}

\begin{figure}[ht]
	\begin{center}
		\centerline{\includegraphics[width=\columnwidth]{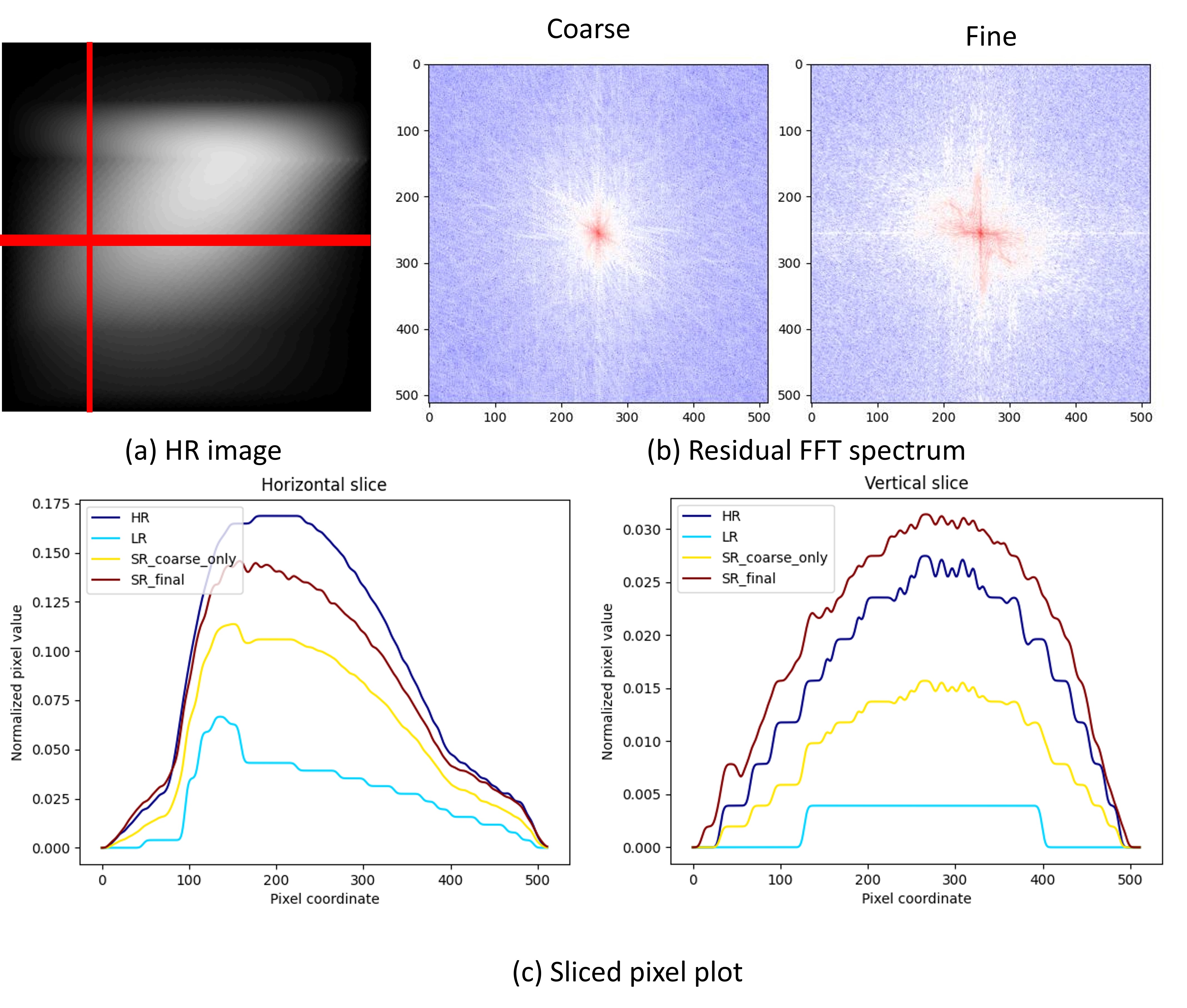}}
        \caption{\small \textbf{Visualization of multiscale feature sampling for 16-times super-resolution.} Given one 16-times super-resolved example ((a) shows the ground truth data), we visualize the frequency map of coarse- and fine-scale reconstructed data in (b), and the plots of sliced pixels in~(c). The slicing position is shown with the red lines. }
		\label{fig:ms}
	\end{center}
\end{figure}

\begin{table}[t]\centering
\caption{\textbf{Ablation study on the multiscale feature sampling and SCS loss.} We report the MSE ($10^{-3}$), MAE ($10^{-3}$) and SSIM results on the test datasets of 4-times and 8-times upsampling.} 
\renewcommand\arraystretch{1.4}
\resizebox{\columnwidth}{!}{
\begin{tabular}{c|cccccc}
\toprule
\multirow{2}{*}{Multiscale} & \multicolumn{3}{c}{4x} & \multicolumn{3}{c}{8x} \\
 & MSE$\downarrow$ & SSIM$\uparrow$ & MAE$\downarrow$ & MSE$\downarrow$ & SSIM$\uparrow$ & MAE$\downarrow$ \\ \midrule
LIIF & 0.874 & 0.857 & 22.16 & 0.952 & 0.922 & 18.86 \\
Single-scale+$\mathrm{L}_1$ & 0.867 & 0.858 & 21.46 & 0.946 & 0.922 & 18.12 \\
Single-scale+$\mathrm{L}_1$+$\mathrm{L}_\mathrm{scs}$ & 0.857 & 0.860 & 21.15 & 0.942 & 0.923 & 17.86 \\
Multiscale+$\mathrm{L}_1$ & 0.848 & 0.859 & 20.86 & 0.907 & 0.920 & 17.45 \\
\cellcolor{mistyrose}{Multiscale+$\mathrm{L}_1$+$\mathrm{L}_\mathrm{scs}$} & \cellcolor{mistyrose}{0.834} & \cellcolor{mistyrose}{0.863} & \cellcolor{mistyrose}{20.09} & \cellcolor{mistyrose}{0.885} & \cellcolor{mistyrose}{0.924} & \cellcolor{mistyrose}{16.93} \\
\bottomrule
\end{tabular}%
}
\label{tab:module}
\end{table}

\begin{figure}[ht]
	\begin{center}
		\centerline{\includegraphics[width=\columnwidth]{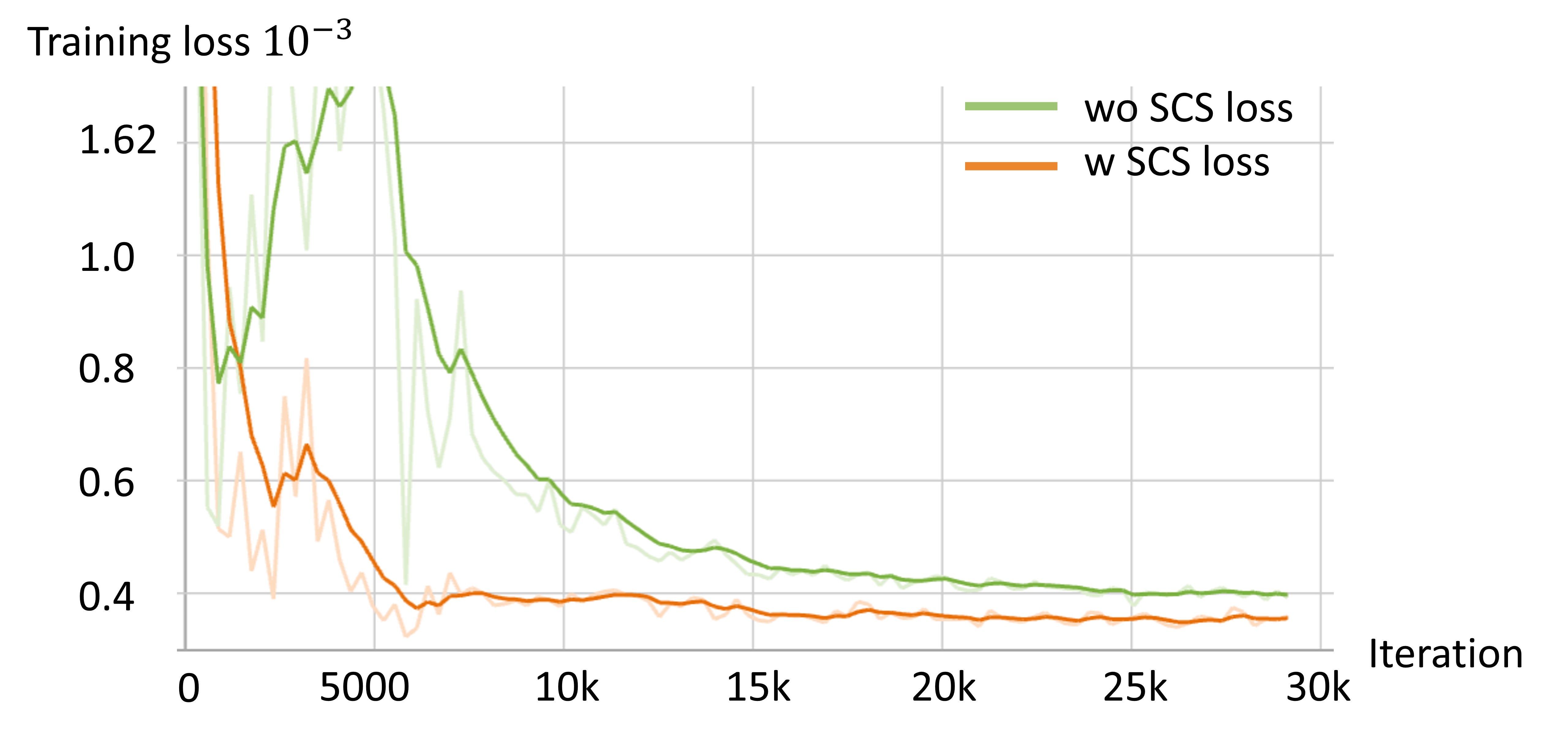}}
        \caption{\small \textbf{Validation loss changes over training iterations.} To demonstrate the effect of the SCS loss for the model optimization, we visualize the validation loss changes over iterations.}
		\label{fig:loss_curve}
	\end{center}
\end{figure}

\noindent \textbf{Multiscale feature sampling and loss optimization.} We conduct an ablation study to show the effect of multiscale feature sampling and the SCS loss for super-resolution. In Table~\ref{tab:module}, we mark \textit{single-scale} as the model that only uses coarse-level features for reconstruction (first equation in~\eqref{eq:multi}), and \textit{multiscale} as the model that uses both coarse- and fine-level features for reconstruction (first and second equations in~\eqref{eq:multi}). LIIF is the baseline for comparison and we can see that using both multiscale feature sampling and the novel SCS loss can improve the reconstruction quality in terms of MSE and SSIM. To further show the effect of using SCS loss, we also plot the validation loss curve in Figure~\ref{fig:loss_curve}. We can see that using the SCS loss, the validation loss is lower and it also converges faster than without using the SCS loss.

\noindent To demonstrate the efficiency of multiscale feature sampling, we show one example of 16-times super-resolution in Figure~\ref{fig:ms}. In (b), we apply the fast Fourier transform (FFT) to the super-resolved data and keep the frequency energy spectrum, and then we shift the zero-frequency component to the center of the array. We can observe that the coarse-scale reconstructed image focuses more on the low-frequency estimation, while the fine-scale reconstruction expands the energy to the wider range of other frequency bands, which indicates that it learns to recover higher-frequency information. In Figure~\ref{fig:ms} (c), we plot the pixel values in the sliced regions (marked by the red lines in (a)). We use the solid dark brown line to represent the results obtained by our final model (combined coarse- and fine-scale feature sampling), and the solid yellow line to represent the results obtained only by coarse-scale feature sampling. We can see that the final result is closer to the ground truth which aligns with the observation in (b) that more details are recovered by the multiscale feature reconstruction. The alignment is still not fully optimal but recall that we are dealing with a complicated setting far beyond the training data. Including further coefficient maps into the training process, can further improve the results.

\begin{figure}[ht]
	\begin{center}
		\centerline{\includegraphics[width=\columnwidth]{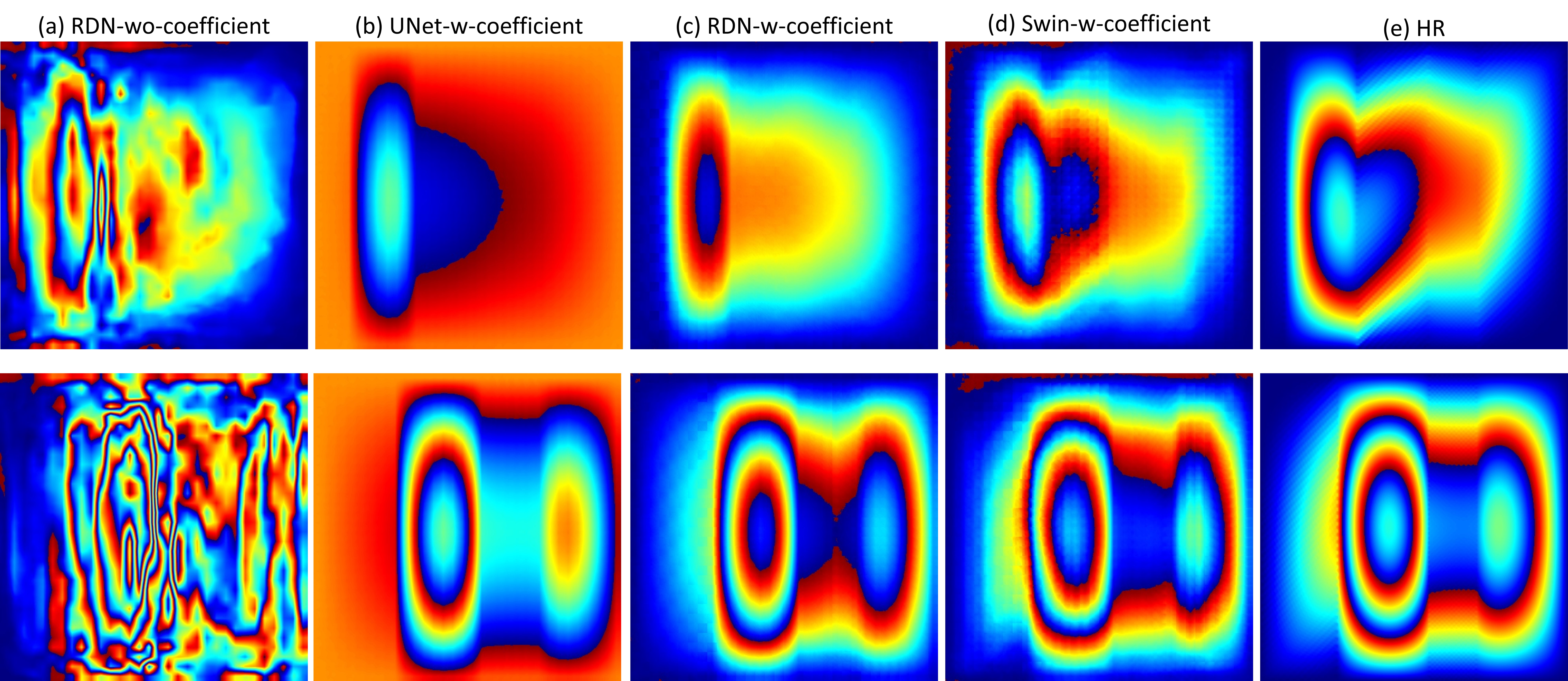}}
        \caption{\small \textbf{Visualization of conditional super-resolution for 32-times super-resolution.} To demonstrate the effect of the coefficient map and different encoder structures, we visualize the super-resolved finite element data. For better comparison, the values are multiplied by a factor 8 and taken modulo 255 as before.}
		\label{fig:condition}
	\end{center}
\end{figure}

\noindent \textbf{Conditional super-resolution.} The proposed \net\ approach takes a coefficient map as the conditional input to jointly learn the super-resolved data. To demonstrate its effect, we provide two examples in Figure~\ref{fig:condition}. Comparing the first and the third column, we may conclude that using coefficient maps as conditions for super-resolution can significantly improve the visual similarity. From the second, third, and fourth row, it becomes clear that using a more complex global feature encoder can improve the reconstruction quality. It is worth noting that using swin transformers~\cite{swinir} (cf.~the fourth column) as the encoder would cost much more computation time than RDN~\cite{rdn}, and it also produces higher MSE losses.

\begin{figure*}[ht]
	\begin{center}
		\centerline{\includegraphics[width=\textwidth]{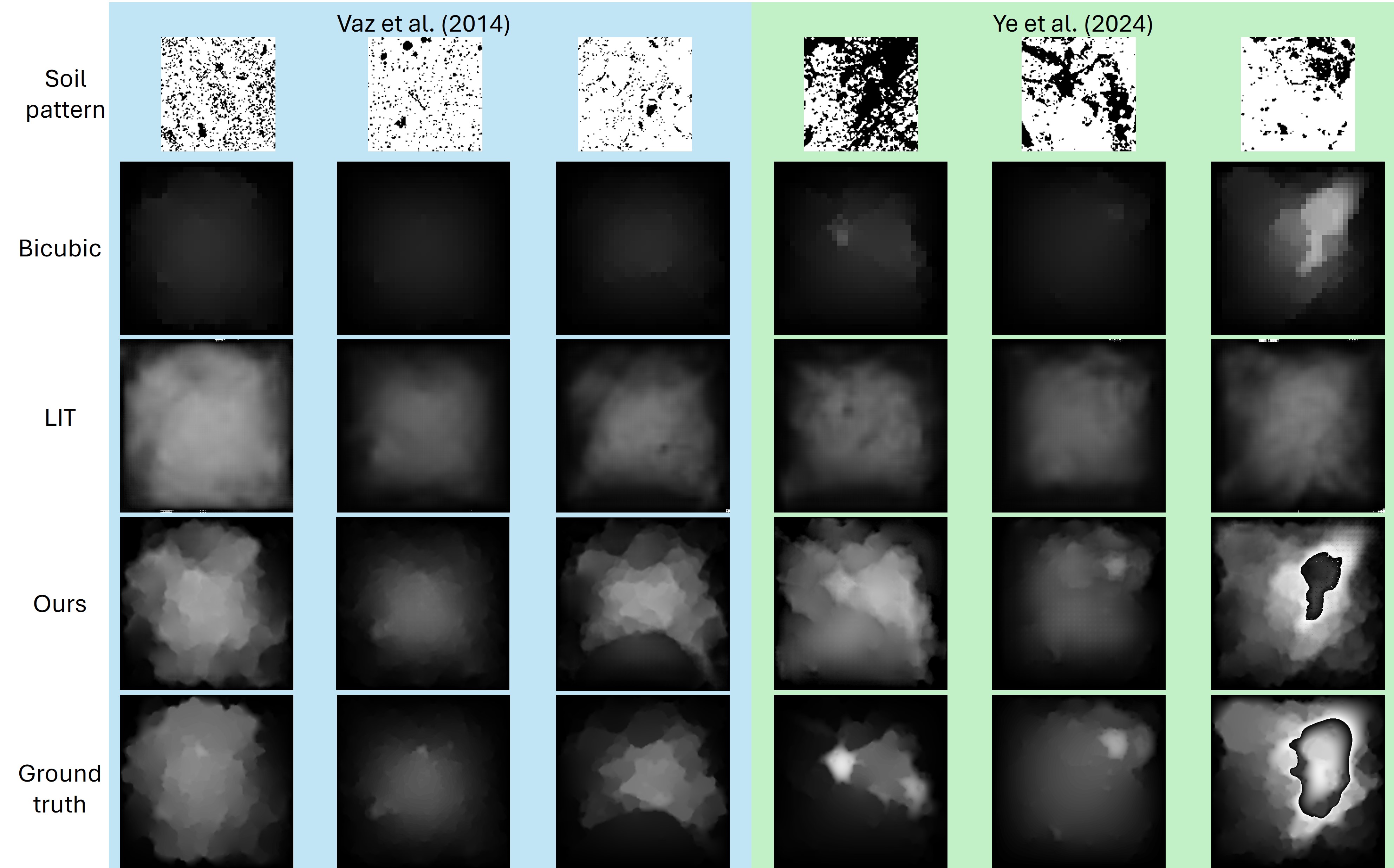}}
        \caption{\small \textbf{Visualization of real-world soil patterns for 16-times super-resolution.} To demonstrate the effect of the proposed method for real-world data, we visualize the super-resolved finite element data. For better comparison, the values are multiplied by a factor of 10 and taken modulo 255 as before.}
		\label{fig:real_world}
	\end{center}
\end{figure*}

\subsection{Real-world applications.} We are interested in applying our proposed method to real-world data to validate its robustness without domain-specific model fine-tuning. In Figure~\ref{fig:real_world}, we show comparisons on two datasets inspired by two research works: 1) Ye et al. (2024)~\cite{YeXZC24} and 2) Vaz et al. (2014)~\cite{VazTLC14}. The first row shows the binary microCT soil patterns used as inputs to our FEM solver to generate paired low-resolution fields ($33\times33$) and high-resolution references (fourth row, $513\times513$). All competing methods are tasked with performing $16$-times super-resolution. Across both datasets, the bicubic baseline produces overly smoothed predictions with little structural detail, while LIT partially restores coarse patterns but still lacks sharp transitions (introducing extra ringing effects) and accurate spatial variability. In contrast, our method (fourth row) recovers markedly clearer boundaries, more faithful spatial gradients, and finer heterogeneous structures that are visually much closer to the ground truth. These improvements are especially evident in regions with sharp contrast or irregular soil clusters, where competing methods tend to oversmooth or distort the shapes. This demonstrates that our approach generalizes well to real-world data and preserves physically meaningful structures even without task-specific retraining.

\section{Conclusion}
\label{Conclusion}
In this paper, we propose a numerical homogenization method based on a continuous super-resolution strategy (\net). Given a coefficient map corresponding to a highly heterogeneous elliptic PDE and a simple coarse finite element approximation, it learns to generate multiscale super-resolved finite element approximation from the continuous feature space. The proposed \net\ strategy can generate results with various resolutions, even beyond the ones used during training. Our proposed Gabor wavelet-based coordinate encoding provides a smoother distribution of grid values and captures the underlying physics more accurately. The multiscale grid feature sampling utilizes the attention-based feature aggregation to effectively learn and predict patterns. {Adding the final stochastic cosine similarity score better aligns the visual appearance to the ground truth than the variant without it.} Remarkably, training only with very specific random coefficients and corresponding finite element solutions already allows for improved (macroscopic) predictions for more general cases (e.g.\ with structured coefficient). Our study thus paves a novel direction to employ neural super-resolution models for numerical homogenization, and we are interested in further exploring and improving ultra-resolution reconstructions. In particular, a thorough choice of more complex training data can already lead to even better generalized predictions.

\section*{Acknowledgments}

Z.~S.~Liu has has been supported by the Research Council of Finland's decisions number 350101 and 359633.

R.~Maier acknowledges support from the German Academic Exchange Service (DAAD), project ID 57711336.

A.~Rupp has been supported by the Research Council of Finland's decisions number 350101, 354489, 359633, 358944, and Business Finland's project number 539/31/2023.

\bibliographystyle{IEEEtran}
\bibliography{main}


\section{Biography Section}
\vspace{-30pt}
\begin{IEEEbiography}[{\includegraphics[width=1in,height=1.25in,clip,keepaspectratio]{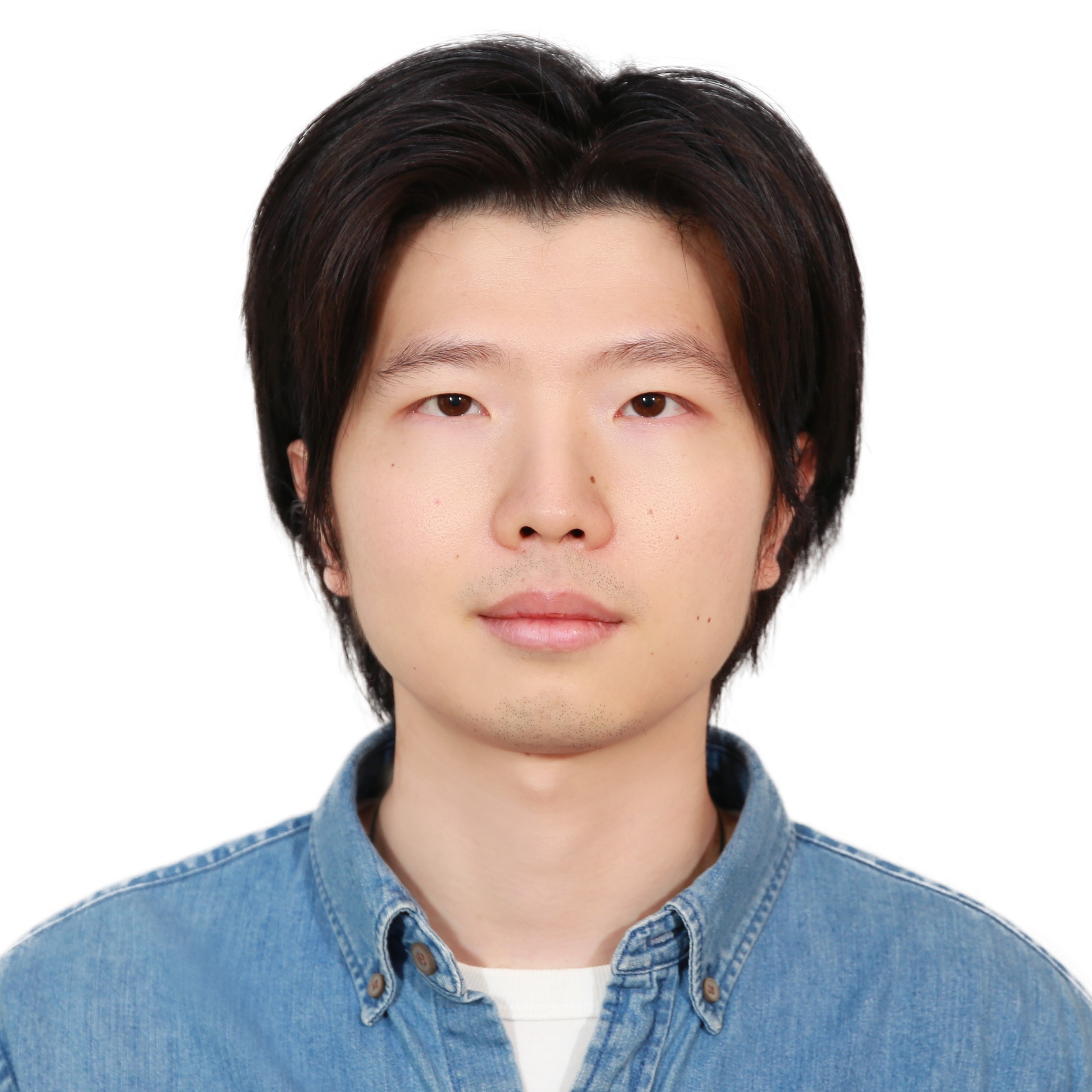}}]{Zhi-Song Liu}
received his Ph.D. (2016-2020) from The Hong Kong Polytechnic University, Artificial Intelligence and Signal Processing Lab. Previously, he was a principal research scientist at DELL EMC. Before that, he was a research fellow (01/2021-12/2021) in CIHE and a post-doc researcher in LIX, Ecole Polytechnique (01/2020-12/2020). His research interests include image and video signal processing, computer vision, deep learning, and 3D data processing. He joined the CVPR Lab at LUT University in 2023 as an associate professor.
\end{IEEEbiography}

\vspace{-15pt}
\begin{IEEEbiography}[{\includegraphics[width=1in,height=1.25in,clip,keepaspectratio]{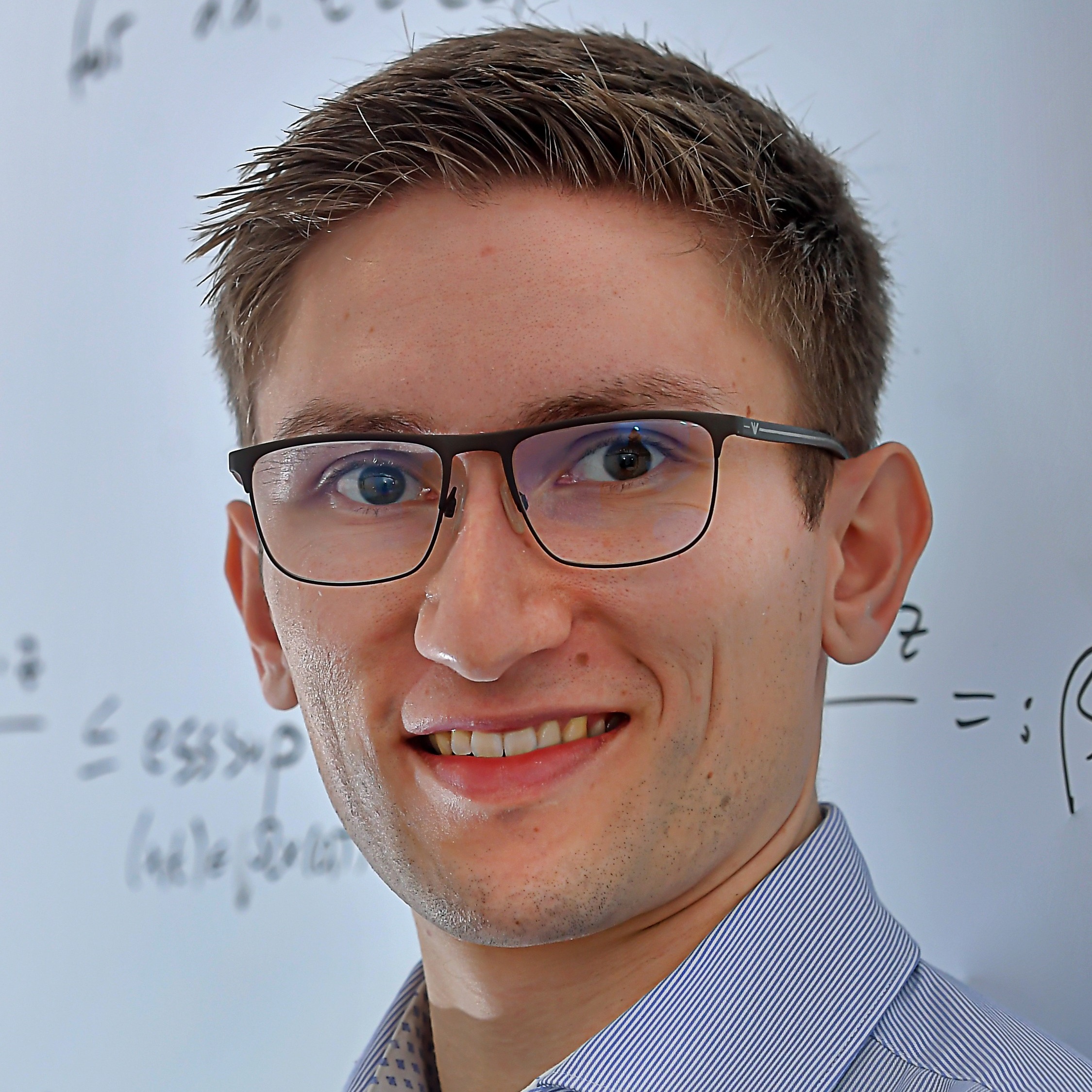}}]{Roland Maier}
studied Mathematics at the University of Bonn (2012-2017) and received his Ph.D. from the University of Augsburg in 2020. After some time as a postdoc at the University of Gothenburg and Chalmers University of Technology, he became a junior professor at University of Jena in 2021. In 2023, he joined Karlsruhe Institute of Technology as a tenure-track professor for numerics of partial differential equations. He mainly works on multi-scale methods and numerical homogenization approaches.   
\end{IEEEbiography}

\vspace{-15pt}
\begin{IEEEbiography}[{\includegraphics[width=1in,height=1.25in,clip,keepaspectratio]{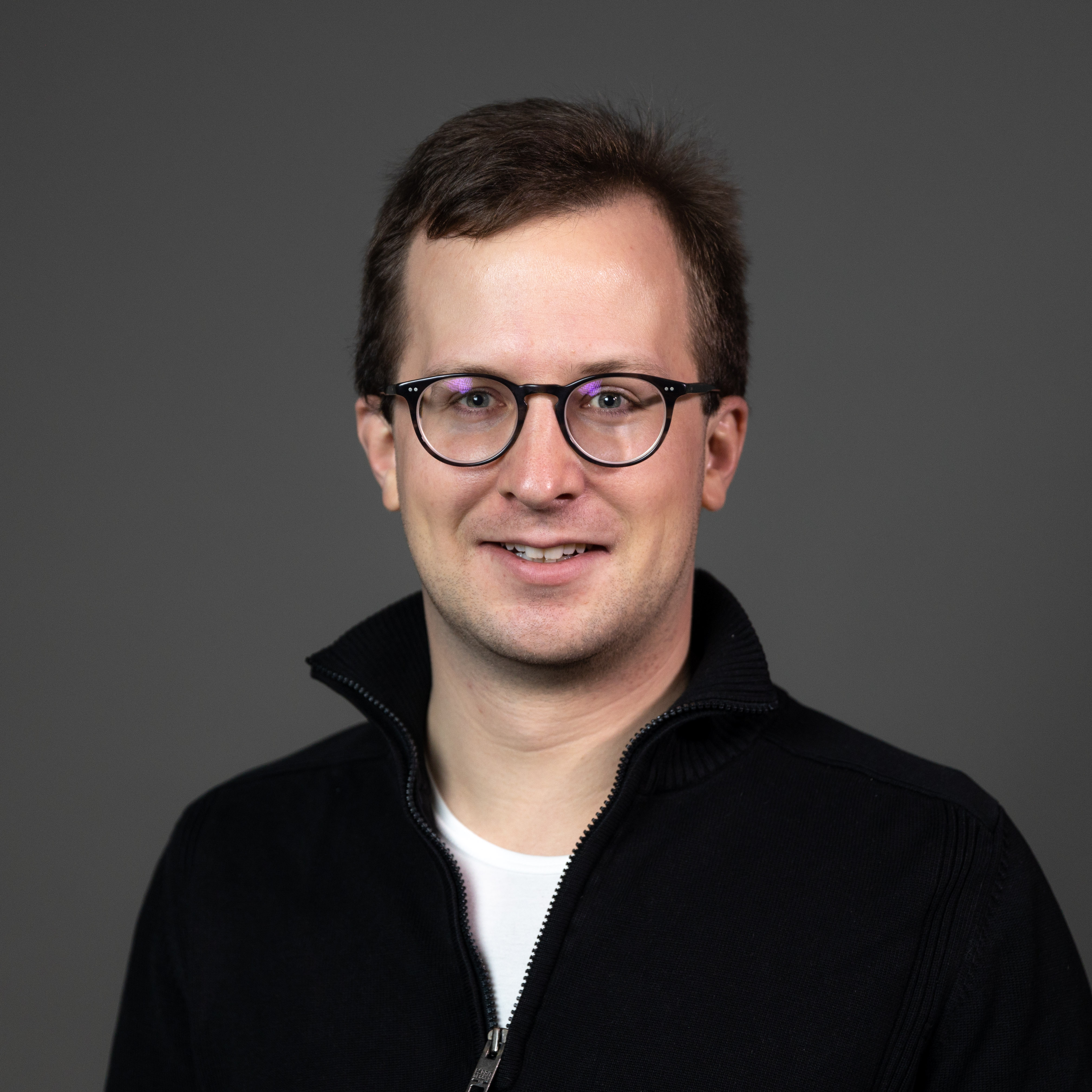}}]{Andreas Rupp}
received his M.Sc.\ and PhD degrees from the FAU Erlangen--Nuremberg in December 2015 and May 2019, respectively. He was a postdoctoral researcher at Heidelberg University till June 2021. Since then, he was a tenure-track assistant professor at LUT University and, since September 2023, an Academy of Finland research fellow. In September 2024, he joined Saarland University as a professor of applied mathematics.
\end{IEEEbiography}

\vfill

\end{document}